\documentclass{siamonline171218}

\usepackage{url}
\usepackage{dsfont}
\usepackage{amsmath, amssymb, color, bm}
\usepackage{graphicx}
\usepackage{float}
\usepackage{mathtools}
\usepackage{enumitem}
\usepackage{indentfirst}
\usepackage{verbatim}
\usepackage{chemfig, chemnum}
\usepackage{tikz}
\usepackage{tikz-cd}
 \usepackage{enumitem}
\usepackage{makecell}
\usepackage{array}
\usepackage{amsopn}

\setlist[enumerate]{leftmargin=.5in}
\setlist[itemize]{leftmargin=.5in}

\newtheorem{example}[theorem]{Example}
\newtheorem{conjecture}[theorem]{Conjecture}

\newcommand{\PP}{\mathbb{P}}
\newcommand{\RR}{\mathbb{R}}

\newcommand{\CC}{\mathbb{C}}
\newcommand{\ZZ}{\mathbb{Z}}

\newcommand{\SH}{\mathrm{SH}}
\newcommand{\Gr}{\mathrm{Gr}}

\newcommand{\Pcal}{\mathcal{P}}

\renewcommand{\tilde}{\widetilde}

\ifpdf
  \DeclareGraphicsExtensions{.eps,.pdf,.png,.jpg}
\else
  \DeclareGraphicsExtensions{.eps}
\fi

\newsiamremark{remark}{Remark}
\newsiamremark{hypothesis}{Hypothesis}
\Crefname{hypothesis}{Hypothesis}{Hypotheses}
\newsiamthm{claim}{Claim}

\headers{Spinor-Helicity Varieties}{Y. El~Maazouz, A. Pfister and B. Sturmfels}

\title{Spinor-Helicity Varieties
\thanks{Submitted to the editors DATE.}}

\author{
        Yassine El Maazouz\thanks{Caltech (\email{maazouz@caltech.edu}).} 
  \and Ana\"elle Pfister\thanks{MPI-MiS Leipzig (\email{anaelle.pfister@mis.mpg.de}).}
  \and Bernd Sturmfels\thanks{MPI-MiS  Leipzig (\email{bernd@mis.mpg.de}).}
}

\makeatletter
\newcommand*{\addFileDependency}[1]{
  \@addtofilelist{#1}
  \IfFileExists{#1}{}{\typeout{No file #1.}}
}
\makeatother

\ifpdf
\hypersetup{
  pdftitle={Spinor-Helicity Varieties},
  pdfauthor={Y. El~Maazouz, A. Pfister and B. Sturmfels}
}
\fi

\begin{document}
\maketitle

\begin{abstract}
\noindent
The spinor-helicity formalism in particle physics gives rise to
natural subvarieties in the product of two Grassmannians.
These include two-step flag varieties for subspaces of
complementary dimension. Taking Hadamard products leads to
Mandelstam varieties.
We study these varieties through the lens of
combinatorics and commutative algebra, 
and we explore their tropicalization,
positive geometry, and scattering correspondence.
 \end{abstract} 



\section{Introduction}

Given two matrices $\lambda $ and $\tilde \lambda$ of
format $k \times n$, suppose
that the $k \times k$ matrix 
$\lambda \cdot \tilde \lambda^T$
has rank at most $r$, for some $0 \leq r \leq k \leq n$.
We wish to express this property in terms of the $k \times k$ minors of the  matrices $\lambda$ and $\tilde \lambda$.
This situation arises in the study of {\em scattering amplitudes}
in quantum field theory \cite{BHPZ}. The special case when
$k = 2$ and $r=0$ is known as {\em spinor-helicity formalism};
for textbook basics see
\cite[Section~1.8]{BHPZ} and
\cite[Section~2.2]{EH}.
In physics, it is customary to write $\langle i j  \rangle$ for the
$2 \times 2$ minors of $\lambda$ 
and $[ i j ]$ for the $2 \times 2$ minors of $\tilde \lambda$,
where $1 \leq i < j \leq n$, and these minors satisfy the
{\em momentum conservation} relations.

\begin{example}[$k=2,n=5,r=0$] \label{ex:twofivezero}
We consider the two skew-symmetric $5 \times 5$ matrices
\[
\!\!  P \, = \, 
\begin{pmatrix}
\,\,0 & \langle 1 2 \rangle & \langle 1 3 \rangle & \langle 1 4 \rangle &\langle 1 5 \rangle \\
\! - \langle 1 2 \rangle & 0 & \langle 2 3 \rangle & \langle 2 4 \rangle &\langle 2 5 \rangle \\
\! - \langle 1 3 \rangle &\!\! - \langle 2 3 \rangle & 0 & \langle 3 4 \rangle &\langle 3 5 \rangle \\
\! - \langle 1 4 \rangle &\!\! - \langle 2 4 \rangle &\!\! - \langle 3 4 \rangle & 0 & \langle 4 5 \rangle \\
\! - \langle 1 5 \rangle &\!\! - \langle 2 5 \rangle &\!\! - \langle 3 5 \rangle &\!\! - \langle 4 5 \rangle & 0 \\
\end{pmatrix} 
\quad {\rm and} \quad Q \, = \,
 \begin{pmatrix}
\,\,0 & [ 1 2 ] & [ 1 3 ] & [ 1 4 ] &[ 1 5 ] \\
\! - [ 1 2 ] & 0 & [ 2 3 ] & [ 2 4 ] &[ 2 5 ] \\
\! - [ 1 3 ] &\!\! - [ 2 3 ] & 0 & [ 3 4 ] &[ 3 5 ] \\
\! - [ 1 4 ] &\!\! - [ 2 4 ] &\!\! - [ 3 4 ] & 0 & [ 4 5 ] \\
\! - [ 1 5 ] &\!\! - [ 2 5 ] &\!\! - [ 3 5 ] &\!\! - [ 4 5 ] & 0 \\
\end{pmatrix}\!.
\]
These matrices have rank two, meaning that the $4 \times 4$ Pfaffians vanish for both matrices:
\begin{equation}
\label{eq:pluecker} \!\!\!
 \langle i j \rangle \langle kl \rangle -
 \langle i k \rangle \langle j l \rangle +
 \langle i l \rangle \langle j k\rangle \,\, = \,\,
 [ i j ] [ kl ] -
 [ i k ] [ j l ] +
 [ i l ] [ j k] \, = \,
 0 \,\,\, \,\hbox{for} \,\,\, 1 \!\leq\! i \!<\! j \!<\! k \!<\! l \!\leq\! 5.
\end{equation}
These {\em quadratic Pl\"ucker relations} are known as {\em Schouten identities} in physics
\cite[eqn (1.116)]{BHPZ}.
Momentum conservation \cite[eqn (1.117)]{BHPZ} stipulates that the product
$P \cdot Q^T$ is the zero matrix:
\begin{equation}
\label{eq:momentum}
\langle i 1 \rangle [1 j] +
\langle i 2 \rangle [2 j] +
\langle i 3 \rangle [3 j] + 
\langle i 4 \rangle [4 j] + 
\langle i 5 \rangle [5 j] \, = \, 0
\quad \hbox{for} \,\, \,1 \leq i,j \leq 5.
\end{equation}
In total, we have a system of $5+5+25 = 35$ quadratic equations
in $\binom{5}{2} + \binom{5}{2} = 20$ unknowns.
The equations (\ref{eq:pluecker}) define a product of
two Grassmannians ${\rm Gr}(2,5) \times {\rm Gr}(2,5) \subset
\PP^9 \times \PP^9$. Inside this product, the bilinear equations (\ref{eq:momentum}) 
cut out a variety
of dimension $8$. This is our spinor-helicity variety, denoted ${\rm SH}(2,5,0)$.
Its bidegree in $\PP^9 \times \PP^9$ is the cohomology class
\begin{equation}
\label{eq:cohom250}  5 s^3 t^7 + 10 s^4 t^6 + 12 s^5 t^5 + 10 s^6 t^4 + 5 s^7 t^3 
\,\, \in \,\, H^*(\PP^9 \times \PP^9, \ZZ). 
\end{equation}
The $35$ quadrics (\ref{eq:pluecker}) and (\ref{eq:momentum})
generate the prime ideal of ${\rm SH}(2,5,0)$. This ideal 
coincides with the
 ideal of the Grassmannian ${\rm Gr}(3,6)$,
which has codimension $10$ and degree $42$ in~$\PP^{19}$.
This identification was observed by Bossinger, Drummond and  Glew \cite{BDG},
who used the term {\em massless scattering ideal} 
and the notation $I_{\rm 5pt}$ for the ideal of ${\rm SH}(2,5,0)$.
 In \cite[Section 6.2]{BDG} they derive the tropicalization of
  ${\rm SH}(2,5,0)$ from that of ${\rm Gr}(3,6)$; see
\cite[Example 4.4.10]{MStrop}.
\end{example}

The Grassmannian ${\rm Gr}(k,n)$ is the subvariety of  $\PP^{\binom{n}{k}-1}$
defined by the Pl\"ucker equations. Points in ${\rm Gr}(k,n)$ represent
$k$-dimensional linear subspaces in~$\CC^n$.
See e.g.~\cite[Chapter~5]{MS}.
Every point in  ${\rm Gr}(k,n)$ is the row space of a
matrix $\lambda \in \CC^{k \times n}$ of rank $k$.
We consider the~set
\begin{equation}
\label{eq:lala1} \bigl\{ \,(\lambda, \tilde \lambda) \,\,:\,\,
\lambda , \tilde \lambda \in \CC^{k \times n} \,,\,\,
{\rm rank}(\lambda) = {\rm rank}(\tilde \lambda) = k,\,\,
{\rm and} \,\,\,
{\rm rank}\bigl( \lambda \cdot \tilde \lambda^T \bigr) \leq r \,\bigr\}.
\end{equation}
We define the {\em spinor-helicity variety} to be
the image of (\ref{eq:lala1}) in the product of  Grassmannians:
\begin{equation}
\label{eq:lala2}
{\rm SH}(k,n,r) \,\,\subset\,\,
{\rm Gr}(k,n) \times {\rm Gr}(k,n) \,\,\subset \,\,
\PP^{\binom{n}{k}-1} \times \PP^{\binom{n}{k}-1}. 
\end{equation}

In what follows, we study the algebra, combinatorics and geometry
of the inclusion (\ref{eq:lala2}). 
Our presentation aims to be accessible to a wide range of readers, not
just from mathematics, but also from physics.
The  prerequisites  are at the level of the
textbooks \cite{MStrop, MS, CCA, AIT}.

The motivation for this project arose from our desire
to understand the spinor-helicity formalism in physics.
The variety ${\rm SH}(2,n,0)$ is widely used
for scattering amplitudes \cite{BHPZ, EH}.
Cachazo, Early, Guevara and Mizera \cite[Section 5.1]{CEGM}
proposed the variety ${\rm SH}(k,n,k-2)$ as a model to
encode kinematic data for particle scattering.
Scattering amplitudes in the CEGM model are computed
by integrating over the moduli space $X(k,n)$ of
$n$ points in~$\PP^{k-1}$.
The articles \cite{ABF, CE, ST} studied the
scattering potential on $X(k,n)$. 
The nonlinear structure of the kinematic data
was highlighted in  Lam's  lectures \cite{Lam2}.
We here examine this in detail.

 The kinematic data  are summarized in the
Mandelstam invariants (\ref{eq:mandelvar2}).
 In the $k=2$ case
from Example \ref{ex:twofivezero},
the moduli space is
$X(2,n) = \mathcal{M}_{0,n}$,  and
 the {\em Mandelstam invariants}~are
\begin{equation}
\label{eq:mandels} s_{ij} \,\, = \,\, \langle i j \rangle [ i j ]. 
\end{equation}
These quantities play the role of the data in the
log-likelihood interpretation of \cite{ST}.
Thus, from the algebraic statistics
 perspective, our topic here is the
geometry of data space.
The fundamental object which underlies this geometry is
the spinor-helicity variety  ${\rm SH}(k,n,r)$.

Our presentation is useful
for readers from physics because it
offers a systematic theory of kinematic
spaces for scattering amplitudes, going well
beyond the familiar 
Grassmannian. The algebra we present 
differs from combinatorial structures that
are known in the community, like cluster coordinates, 
and it opens new vistas on positivity and
tropicalization. Our varieties (\ref{eq:lala2})
for $r=k-2$ capture the
essence of matrix kinematics in CEGM theory
\cite[Section~5]{CEGM}.

The article is organized as follows.
In Section~\ref{sec2} we present quadratic polynomials
that form a Gr\"obner basis for the
 prime ideal of ${\rm SH}(k,n,r)$.
 The underlying toric degeneration is represented
 by a poset constructed from two copies of Young's lattice
 for ${\rm Gr}(k,n)$.
 The special case $r=0$ is understood by identifying
${\rm SH}(k,n,0)$ with the
two-step flag variety ${\rm Fl}(k,n-k; \CC^n)$.
 In Section~\ref{sec3}
we express the momentum conservation equations by
a matrix product $P Q^T$ which generalizes (\ref{eq:momentum}), and we show that these equations generate the prime ideal of (\ref{eq:lala2}). Theorem  \ref{thm:khovanskii}
features a Khovanskii basis for the coordinate ring of ${\rm SH}(k,n,r)$.

Section~\ref{sec4} investigates the polynomial relations among the Mandelstam invariants.
These relations define the {\em Mandelstam variety} ${\rm M}(k,n,r)$ 
 in the kinematic subspace of $\PP^{\binom{n}{k}-1}$. 
 We study both parametric and implicit representations of
  this variety. The   generators of its prime ideal for the case $k=2$
  are presented in Theorem \ref{thm:prime2n}.
   We note in Proposition \ref{prop:deletion} 
  that ${\rm M}(k,n,k)$ is the Hadamard
  product (see \cite{BC}) of the Grassmannian ${\rm Gr}(k,n )$ with itself.

In Section \ref{sec5} we turn to positive  geometry and tropical geometry. 
We introduce the positive parts of ${\rm SH}(k,n,r)$
and ${\rm M}(k,n,r)$, we discuss their boundaries,
and we compute some associated tropical varieties.
For $r=0$, these structures arise 
from the flag variety.

In Section~\ref{sec6} we study the scattering correspondence.
This is a variety in the product~space
$$ {\rm M}(k,n,r) \,\times \, X(k,n) , $$
where $X(k,n)$ is the  moduli space for $n$ points in $\PP^{k-1}$. 
It parametrizes pairs of Mandelstam invariants and
solutions to their scattering equations; see \cite[eqn (0.2)]{Lam2}.
This mirrors the likelihood correspondence
in statistics \cite[Definition 1.5]{LG}. Building on \cite[Section 4.3]{Lam2}, we offer a mathematical perspective on results from the physics literature, mostly for $k=2$.

This article is accompanied by software and data. These materials
are made available in the {\tt MathRepo} collection at MPI-MiS via \url{https://mathrepo.mis.mpg.de/SpinorHelicity}.

\section{Two Grassmannians and their Posets}
\label{sec2}

We fix two copies of the Grassmannian 
${\rm Gr}(k,n)$ embedded in $\PP^{\binom{n}{k}-1}$.
The first Grassmannian has Pl\"ucker coordinates
$\langle i_1  i_2 \ldots i_k \rangle$, representing
maximal minors of $\lambda$.
The second one has Pl\"ucker coordinates
$[ i_1  i_2 \ldots i_k ]$, representing
maximal minors of $\tilde \lambda$.
These expressions are antisymmetric, so we usually assume
$1 \leq i_1 < i_2 < \cdots<  i_k \leq n$.
 For instance, for $k = 3$, 
 \begin{equation}
 \label{eq:sorttt}
  \begin{matrix} &
 \langle 123 \rangle 
 = - \langle 132 \rangle 
 = - \langle 213 \rangle 
 =\, \langle 231 \rangle
 = \, \langle 312 \rangle
   \, = - \langle 321 \rangle \\ {\rm and} \quad &
\, [ 123 ]  \,
 = - [ 132 ] \,
 =\, - [ 213 ] \,
 \,=\, [ 231 ]\,
 = \, \,[ 312 ]\,
    = - [ 321 ] .
    \end{matrix}
\end{equation}
Their relations are given by two copies of the Pl\"ucker ideal, denoted
$J_{k,n}$ and ${\tilde J}_{k,n}$.
To describe this ideal, we introduce Young's lattice $Y_{k,n}$.
This is the partially ordered set (poset) whose elements are 
the $\binom{n}{k}$ Pl\"ucker coordinates. The order relation in $Y_{k,n}$ is
defined by
 \begin{equation}\label{eq:YoungLattice_Order1}
 \langle i_1 i_2 \cdots i_k \rangle \, \leq \, 
  \langle j_1 j_2 \cdots j_k \rangle \quad : \Longleftrightarrow \quad
 i_1 \leq j_1 \,\,{\rm and} \,\,
  i_2 \leq j_2 \,\,{\rm and} \,\,
\cdots
  \,\,{\rm and} \,\,
 i_k \leq j_k.    
 \end{equation}\label{eq:YoungLattice_Order2}
 Let $\tilde Y_{k,n}$
 be a second copy of Young's poset, but now with the 
 order relation reversed:
 \begin{equation}
 [ i_1 i_2 \cdots i_k ]\, \leq \, 
 [j_1 j_2 \cdots j_k] \quad : \Longleftrightarrow \quad
 i_1 \geq j_1 \,\,{\rm and} \,\,
  i_2 \geq j_2 \,\,{\rm and} \,\,
\cdots
  \,\,{\rm and} \,\,
 i_k \geq j_k.    
 \end{equation}
 The following result on the ideal $J_{k,n}$ of the Grassmannian
 is well-known (see~\cite[\S3.1]{AIT}).
  
 \begin{proposition} \label{prop:grass}
 The prime ideal $J_{k,n}$ is generated by  the quadratic Pl\"ucker relations
\begin{equation}
\label{eq:pluecker1} \sum_{s=0}^k \,(-1)^s \cdot 
\langle \,i_1 \,i_2 \,\cdots\, i_{k-1} \,j_s \, \rangle \cdot
\langle j_0 j_1 \cdots j_{s-1} j_{s+1} \cdots j_k \rangle . 
\end{equation}
These quadrics are a Gr\"obner basis for the
reverse lexicographic term order
given by any~linear extension of $Y_{k,n}$.
The initial ideal of $J_{k,n}$ is generated by the incomparable pairs in~$Y_{k,n}$.
 \end{proposition}

The key point of this result is that every incomparable pair lifts
to a quadric in $J_{k,n}$.

\begin{example}\label{ex:foureight}[$k=4,n=8$]
The elements $\langle 1278 \rangle$ and $\langle 3456 \rangle$ are incomparable in $Y_{4,8}$. 
The Pl\"ucker relation (\ref{eq:pluecker1}) for the triple
$i_1i_2i_3 = 127$ and the quintuple $j_0j_1 j_2 j_3 j_4 = 34568$~is
\[ 
{\bf \langle 1278 \rangle\langle3456 \rangle}
\,+\,\langle 1267 \rangle\langle3458 \rangle
- \langle1257 \rangle\langle3468 \rangle
+\langle1247 \rangle\langle3568 \rangle
- \langle1237 \rangle\langle4568 \rangle.
\]
The monomials are listed in the reverse lexicographic term order given by any linear extension of $Y_{k,n}$.
The initial monomial is the prescribed incomparable pair. 
Our quadric is not in the reduced Gr\"obner basis since it~has incomparable trailing terms.
The corresponding element in the reduced Gr\"obner basis  equals
 \begin{small} $$ 
 \begin{matrix}
 {\bf
  \langle1278 \rangle \langle3456 \rangle
- \langle1256 \rangle \langle3478 \rangle}
+ \langle1246 \rangle \langle3578 \rangle
- \langle1245 \rangle \langle3678 \rangle\quad\qquad \\ \qquad \qquad \qquad
- \,\langle1236 \rangle \langle4578 \rangle
\,+\, \langle1235 \rangle \langle4678 \rangle
\,-\, \langle1234 \rangle \langle5678 \rangle .
\end{matrix}
$$ \end{small}
The quadric above has the virtue that its initial binomial
$  \langle1278 \rangle \langle3456 \rangle - \langle1256 \rangle \langle3478 \rangle$ is consistent with
the toric degeneration of the Grassmannian  ${\rm Gr}(4,8)$
given by Young's lattice~$Y_{4,8}$.
Indeed, for the incomparable pair
$\langle1278 \rangle \langle3456 \rangle$, the meet is $\langle3478 \rangle$ and the join is $\langle1256 \rangle$.
Algebraically, this is the Khovanskii basis (or SAGBI basis) structure
 in  \cite[Theorem 3.2.9]{AIT}.
\end{example}

\begin{corollary}
The number of  generators for $J_{k,n}$, or 
 of incomparable pairs in $Y_{k,n}$, 
equals
\begin{equation}
\label{eq:hookcontent}
\frac{1}{2} \left[\binom{n}{k}+1 \right]\binom{n}{k} \,\, - \,\,
\frac{(n+1)(n-k+1)}{k+1} \prod_{i=0}^{k-2} \frac{(n-i)^2}{(k-i)^2}.
\end{equation}
\end{corollary}

\begin{proof}
The first term is the number of all quadratic monomials in the
$\binom{n}{k}$ variables $\langle i_1 i_2 \cdots i_k \rangle$.
From this we subtract the number of standard monomials,
which is the number of semi-standard Young tableaux
of shape $k \times 2$ with fillings in $\{1,2,\ldots,n\}$.
That number is given by the {\em hook-content formula}
from combinatorics, which we made explicit in (\ref{eq:hookcontent}).
\end{proof}

We now turn to the spinor-helicity variety.
Rephrasing the definition in~(\ref{eq:lala2}),
this~is
\begin{equation} \label{eq:rephrasing}
    \SH(k,n,r) \,\,= \,\,\bigl\{ \,(V,W) \in \Gr(k,n) \times \Gr(k,n) \,\,\colon \dim(V \cap W^\perp ) \geq k-r
\,     \bigr\},
     \end{equation}
where $W^\perp$ is the space orthogonal to $W$ with respect to the 
standard inner product on $\CC^n$.
In the setting of the Introduction,
$V$ and $W$ are the row spaces of $\lambda$ and $\tilde \lambda$
respectively.

\begin{remark}[Involution] \label{rmk:involution}
There is a canonical involution on $\SH(k,n,r)$, defined by
swapping the subspaces $V$ and $W$.
This  interchanges the coordinates
$\langle i_1 i_2 \cdots i_k \rangle$
and $[i_1 i_2 \cdots i_k ]$.
\end{remark}

\begin{proposition} \label{prop:geoflag}
Fix integers $k,n,r$ such that $ 0 \leq r \leq k$ and $2k \leq r+n$.
The spinor-helicity variety $\, \SH(k,n,r)\,$
 is non-empty and irreducible in
  $  \PP^{\binom{n}{k}-1} \times \PP^{\binom{n}{k}-1}$. Its dimension~equals
\begin{equation}
\label{eq:dimformula}
    \dim(\SH(k,n,r)) \,\,=\,\, 2k(n-k) - (k-r)^2.
\end{equation}
If $r=0$ then it is linearly isomorphic to the
two-step flag variety $\,{\rm Fl}(k,n-k;\CC^n)$.
\end{proposition}

\begin{proof}
Our hypothesis that
  $k = {\rm dim}(V)$ and
 $ n-k = {\rm dim}(W^\perp)$ are at least $ k-r$ is 
necessary for $\SH(k,n,r)$ to be non-empty.
We assume $k-r \geq 0$ to rule out trivial cases.
The projection of $\SH(k,n,r)$ onto the first factor
equals ${\rm Gr}(k,n)$, which is irreducible of dimension $k(n-k)$.
The fibers are subvarieties in the second factor ${\rm Gr}(k,n)$.
In local affine coordinates, the fiber over $V$ consists of all $k \times (n-k)$ matrices of rank at most $r$. 
This is an irreducible variety of codimension $(k-r)^2$.
Hence, each fiber is irreducible of the same dimension $k(n-k)-(k-r)^2$.
From this we conclude that $\SH(k,n,r)$ is irreducible of dimension $2k(n-k) - (k-r)^2$.
 
Fix $r = 0$ and recall that $k \leq n-k$.
By passing from $W$ to $W^\perp$, we can view ${\rm SH}(k,n,0)$
as a subvariety in ${\rm Gr}(k,n) \times {\rm Gr}(n-k,n)$.
Its points are pairs $(V,W^\perp)$ of linear subspaces 
in $\CC^n$ such that $V \subseteq W^\perp$.
In other words, its points are two-step flags 
 in $\CC^n$. Hence
${\rm SH}(k,n,0)$ coincides with  the flag variety ${\rm Fl}(k,n-k;\CC^n)$
in its Pl\"ucker embedding in 
$\,   \PP^{\binom{n}{k}-1} \times \,\PP^{\binom{n}{k}-1}$.
\end{proof}

Let $R = \mathbb{C}\bigl[\langle i_1 \dots i_k \rangle , [j_1 \dots j_k]\bigr]$ 
be the polynomial ring in the $2 \binom{n}{k}$ bracket variables.
Let $S = \mathbb{C}[\bf{x}]$ be the polynomial ring 
in the entries of an $(n-k+r) \times n$ 
matrix ${\bf x} = (x_{ij})$. 
We write $\phi_{k,n,r}: R \rightarrow S$ for the
homomorphism which maps
$\langle I \rangle = \langle i_1 i_2 \cdots i_k \rangle$
to the $k \times k$ minor of ${\bf x}$ in the rows
$1,2,\ldots,k$ and  columns $I$,
and  which maps
$[J] = [j_1 j_2 \cdots j_k]$
to $(-1)^{j_1+j_2+\cdots + j_k}$ times the $(n\!-\!k) \times (n\!-\!k)$ minor of ${\bf x}$ in the rows
$r+1,\ldots,n-k+r$ and  columns $[n] \backslash J$.
Note that all such minors involve the $k-r$ middle rows,
which are indexed by $r+1,\ldots,k$.

\begin{remark}[Parametrization] \label{rmk:kernel}
Let $I_{k,n,r}$ denote the kernel of   the ring map $\phi_{k,n,r} $.
This kernel is a homogeneous prime ideal in $R$, and its
zero set is the spinor-helicity variety ${\rm SH}(k,n,r)$.
 Hence, the subalgebra
 $\phi_{k,n,r}(R)$ of $S = \CC[{\bf x}]$ is isomorphic to the 
 coordinate ring of ${\rm SH}(k,n,r)$. Indeed,
following (\ref{eq:rephrasing}), the minors
$\phi_{k,n,r}(\langle I \rangle) $ are the Pl\"ucker coordinates
for the space $V$, while the signed minors
$\phi_{k,n,r}([J])$ are  the Pl\"ucker coordinates
for  $W^\perp$. The 
condition ${\rm dim}(V \cap W^\perp) \geq k-r$
is encoded by the overlap in the $k-r$ middle rows of ${\bf x}$.
\end{remark}

We will describe a Gr\"obner basis of quadrics for $I_{k,n,r}$.
The initial monomials  admit
a combinatorial description which extends that for
Grassmannians given in
Proposition \ref{prop:grass}.
We start out with our two copies of
Young's lattice, $Y_{k,n}$ and $\tilde Y_{k,n}$.
We define a new poset $\mathcal{P}_{k,n,r}$
as follows. As a set, $\mathcal{P}_{k,n,r}$
is the disjoint union of $Y_{k,n}$ and $\tilde Y_{k,n}$.
All order relations in $Y_{k,n}$ and $\tilde Y_{k,n}$
remain order relations in $\mathcal{P}_{k,n,r}$.
In addition, there are $\binom{2k-2r}{k-r}$ covering relations
\begin{equation}
\label{eq:coveringrel}
 [1 2 \cdots r \,i_{r+1}\cdots i_k] \, \leq \,\langle 12 \cdots r \,j_{r+1} \cdots j_k \rangle, 
 \end{equation}
one for each ordered set partition
$\, \{r+1,r+2,\ldots,2k -r\} = \{i_{r+1} ,\ldots, i_k\} \sqcup \{j_{r+1} ,\ldots, j_k \} $.
The poset $\mathcal{P}_{k,n,r}$ is the transitive closure of these relations.
Note that $\mathcal{P}_{k,n,r}$ is a graded poset, with unique
minimal element  $[n{-}k{+}1 \,\cdots \,n{-}1 \, n]$
and unique maximal element
$\langle n{-}k{+}1 \,\cdots \,n{-}1 \, n \rangle$.
The Hasse diagram of $\mathcal{P}_{k,n,r}$ 
is shown in Figure \ref{fig:P_(2,6,0)}
for $k=2, n= 6, r = 0$.

\begin{figure}[!htbp]
    \centering
    \scalebox{0.7}{
     \begin{tikzpicture}

        \draw (2,4) node  {$\langle 56 \rangle$};
        \draw (1,3) node  {$\langle 46 \rangle$};
        \draw (0,2) node  {$\langle 36 \rangle$};
        \draw (-1,1) node {$\langle 26 \rangle$};
        \draw (-2,0) node {$\langle 16 \rangle$};
        
        \draw (1.25,3.25)--(1.75,3.75);
        \draw (0.25,2.25)--(0.75,2.75);
        \draw (-0.75,1.25)--(-0.25,1.75);
        \draw (-1.75,0.25)--(-1.25,0.75);

        \draw (1.25,2.75)--(1.75,2.25);
        \draw (0.25,1.75)--(0.75,1.25);
        \draw (-0.75,0.75)--(-0.25,0.25);
        \draw (-1.75,-0.25)--(-1.25,-0.75);

        \draw (2,2) node {$\langle 45 \rangle$};
        \draw (1,1) node {$\langle 35 \rangle$};
        \draw (0,0) node {$\langle 25 \rangle$};
        \draw (-1,-1) node {$\langle 15 \rangle$};

        \draw (1.25,1.25)--(1.75,1.75);
        \draw (0.25,0.25)--(0.75,0.75);
        \draw (-0.75,-0.75)--(-0.25,-0.25);

        \draw (1.25,0.75)--(1.75,0.25);
        \draw (0.25,-0.25)--(0.75,-0.75);
        \draw (-0.75,-1.25)--(-0.25,-1.75);

        \draw (2,0) node {$\langle 34 \rangle$};
        \draw (1,-1) node {$\langle 24 \rangle$};
        \draw (0,-2) node {$\langle 14 \rangle$};

        \draw (1.25,-0.75)--(1.75,-0.25);
        \draw (0.75,-1.25)--(0.25,-1.75);

        \draw (1.25,-1.25)--(1.75,-1.75);
        \draw (0.25,-2.25)--(0.75,-2.75);

        \draw (2,-2) node {$\langle 23 \rangle$};
        \draw (1,-3) node {$\langle 13 \rangle$};

        \draw (1.75,-2.25)--(1.25,-2.75);

        \draw (1.25,-3.25)--(1.75,-3.75);
         
        \draw (2,-4) node {$\langle 12 \rangle$};

        \begin{scope}[shift={(7,-4)}, rotate=180]

        \draw (2,4) node  {[56]};
        \draw (1,3) node  {[46]};
        \draw (0,2) node  {[36]};
        \draw (-1,1) node {[26]};
        \draw (-2,0) node {[16]};

        \draw (1.25,3.25)--(1.75,3.75);
        \draw (0.25,2.25)--(0.75,2.75);
        \draw (-0.75,1.25)--(-0.25,1.75);
        \draw (-1.75,0.25)--(-1.25,0.75);

        \draw (1.25,2.75)--(1.75,2.25);
        \draw (0.25,1.75)--(0.75,1.25);
        \draw (-0.75,0.75)--(-0.25,0.25);
        \draw (-1.75,-0.25)--(-1.25,-0.75);

        \draw (2,2) node {[45]};
        \draw (1,1) node {[35]};
        \draw (0,0) node {[25]};
        \draw (-1,-1) node {[15]};

        \draw (1.25,1.25)--(1.75,1.75);
        \draw (0.25,0.25)--(0.75,0.75);
        \draw (-0.75,-0.75)--(-0.25,-0.25);

        \draw (1.25,0.75)--(1.75,0.25);
        \draw (0.25,-0.25)--(0.75,-0.75);
        \draw (-0.75,-1.25)--(-0.25,-1.75);

        \draw (2,0) node {[34]};
        \draw (1,-1) node {[24]};
        \draw (0,-2) node {[14]};

        \draw (1.25,-0.75)--(1.75,-0.25);
        \draw (0.75,-1.25)--(0.25,-1.75);

        \draw (1.25,-1.25)--(1.75,-1.75);
        \draw (0.25,-2.25)--(0.75,-2.75);

        \draw (2,-2) node {[23]};
        \draw (1,-3) node {[13]};

        \draw (1.75,-2.25)--(1.25,-2.75);

        \draw (1.25,-3.25)--(1.75,-3.75);
         
        \draw (2,-4) node {[12]};

        \end{scope}


        \draw[red] (2.35,0)--(4.65,0);
        \draw[red] (1.35,-1)--(5.65,-1);

        \draw[red] (2.35,-2) .. controls (4,-1.5) and (5,-1.5) .. (6.65,-2);
        \draw[red] (0.35,-2) .. controls (2,-2.5) and (3,-2.5) .. (4.65,-2);

        \draw[red] (1.35,-3) -- (5.65,-3);
        \draw[red] (2.35,-4) -- (4.65,-4);
        
    \end{tikzpicture}
    }
    \caption{The poset $\Pcal_{2,6,0}$ is created from
        $Y_{2,6}$ and $\tilde Y_{2,6}$ by adding six covering relations.}
    \label{fig:P_(2,6,0)}
\end{figure}

We are now prepared to state our first theorem on the
spinor-helicity variety ${\rm SH}(k,n,r)$.

\begin{theorem} \label{thm:main1}
 The prime ideal $I_{k,n,r}$ is minimally generated by quadratic forms.
 These quadrics are a Gr\"obner basis for the reverse lexicographic term order
given by any~linear extension of $\mathcal{P}_{k,n,r}$.
The initial ideal of $I_{k,n,r}$ is generated by the incomparable pairs in~$\mathcal{P}_{k,n,r}$.
\end{theorem}

\begin{proof} We first assume $r=0$.
By Proposition \ref{prop:geoflag},
our ideal $I_{k,n,0}$ is the ideal of the two-step flag variety
${\rm Fl}(k,n-k; \CC^n)$. The quadratic Gr\"obner basis for
that ideal is derived from the well-known straightening law for
flag varieties. We refer to \cite[Chapter 14]{CCA}
for a textbook exposition.
That exposition emphasizes the case of the
complete flag variety.  This case applies to our situation as follows.
Let $\mathcal{P}$ be the poset 
on all $2^n$ subsets of $\{1,2,\ldots,n\}$ that was introduced in
 \cite[Section 14.2]{CCA}. The restriction of that poset 
 to subsets that have size $k$ or $n-k$ is isomorphic to our
 poset $\mathcal{P}_{k,n,0}$. The poset isomorphism maps
  $(n-k)$-sets to their complements. With this,
  our assertion for $r=0$ follows 
from \cite[Theorem 14.6]{CCA}.

We next present the proof for $r \geq 1$. This will generalize
the known construction we
used for $r=0$. We consider the skew Young diagram
$\lambda/\mu$ where $\lambda = (n-k+r,r)$ and $\mu = (r)$.
A filling of $\lambda/\mu$ with entries in $[n] = \{1,2,\ldots,n\}$
is assumed to have its rows strictly increasing. Hence
there are $\binom{n}{k}^2$ such fillings. A filling
is {\em semi-standard} if the 
$k-r$ non-trivial columns are weakly increasing.
If this is not the case then the filling of $\lambda/\mu$ is called {\em non-standard}.

With these definitions in place, our poset admits the following alternative description:
\begin{equation}
\label{eq:3equivalent}
 \begin{matrix} & \quad \langle i_1 i_2 \,\cdots i_k\, \rangle \,\geq \, [j_1 j_2 \,\cdots \, j_k]  
\qquad \qquad \hbox{holds in $\mathcal{P}_{k,n,r}$,}  \\
\iff & \langle i_{r+1} \!-\! r \,\cdots\, i_k\! - \! r \rangle \,\geq \, [j_{r+1} \!-\! r\, \cdots\, j_k \!-\! r ]
\quad \hbox{holds in $\mathcal{P}_{k-r,n-r,0}$,} \\
\iff & \hbox{the filling of $\lambda/\mu$ with $[n] \backslash \{j_1,\ldots,j_k\}$
and   $\{i_1,\ldots,i_k\}$ is semi-standard.}
\end{matrix} 
\end{equation}
The proof of our theorem is now analogous to that of \cite[Theorem 14.6]{CCA}.
Fix an incomparable pair $\langle I \rangle [J]$ in 
the poset $\mathcal{P}_{k,n,r}$, and consider the
corresponding non-standard skew tableaux
\setcounter{MaxMatrixCols}{11}
\begin{equation}\label{eq:skewSymTableau}
\lambda/\mu \quad  = \qquad
  \begin{bmatrix} &  &  & j_1' & \cdots &  j'_l & \cdots  &j'_{k-r} & j'_{k-r+1} & \cdots & j'_{n-k} \\
               i_1  & \cdots & i_r & i_{r+1} & \cdots & i_{r+l} & \cdots & i_k & & & \\
               \end{bmatrix} .    
\end{equation}
The rows are increasing, and $i_{r+l} < j'_{l}$ is the leftmost violation, and
$\{j'_1,\ldots,j'_{n-k}\} = [n]\backslash J$.

By summing over all permutations $\pi$ of 
$i_1 < \cdots < i_{r+l} < j'_{l} < \cdots < j'_{n-k}$, we obtain
\begin{equation}
\label{eq:pialternating}
 \sum_{\pi} {\rm sign}(\pi) \cdot \langle \pi(I) \rangle \cdot  [\pi([n]\backslash J)]  \quad
 \in \,\, R.
 \end{equation}
This is the analogue to \cite[eqn (14.2)]{CCA}.
The image of (\ref{eq:pialternating}) under $\phi_{k,n,r}$
is an alternating multilinear form in $n-k+r+1$ column vectors
of the matrix ${\bf x}$. The matrix has only $n-k+r$ rows,
so this multilinear form is zero. Therefore (\ref{eq:pialternating}) lies
in $I_{k,n,r} = {\rm kernel}(\phi_{k,n,r})$.
Finally, we note that the initial monomial of (\ref{eq:pialternating}) is
the monomial $\langle I \rangle [J]$ we started out with.

To complete the proof, we need to show that the
semi-standard monomials in $R$ are
linearly independent modulo $I_{k,n,r}$.
The argument for this follows that in the proof of
\cite[Theorem 14.6]{CCA}. Namely, we consider
any monomial in $S$ and we write it as in
\cite[eqn (14.4)]{CCA}. There exists a unique
semi-standard skew tableau whose image
under $\phi_{k,n,r}$ has that initial monomial.
Thus, no cancellation is possible, and this finishes the proof
of Theorem \ref{thm:main1}.
\end{proof}
 
\begin{example}[$k=2,n=6,r=0$]
The ideal $I_{2,6,0}$ is generated by
$15+15+36 = 66$ quadrics which form a Gr\"obner basis.
Their initial monomials are the incomparable pairs
in the poset $\mathcal{P}_{2,6,0}$ which is shown in
Figure \ref{fig:P_(2,6,0)}. 
The $15$ initial monomials from $J_{2,6}$ are the
pairs $\langle ij \rangle \langle kl \rangle $ in $Y_{2,6}$.
The $15$ initial monomials from $\tilde J_{2,6}$ are the
pairs $[ ij ][  kl ]$ in $\tilde Y_{2,6}$.
Finally, there are $36$ mixed initial monomials
$\langle ij \rangle [ kl]$, corresponding to bilinear
generators of $I_{2,6,0}$.

The poset $\mathcal{P}_{2,5,0}$ arises in Figure \ref{fig:P_(2,6,0)} from deleting the 
upper rim
$\langle 1 6 \rangle ,\langle 2 6 \rangle ,\langle 3 6 \rangle ,\langle 4 6 \rangle ,\langle 5 6 \rangle $
and the lower rim $[ 1 6 ] ,[ 2 6 ] ,[ 3 6 ] ,[ 4 6 ] ,[ 5 6 ] $.
Thus $\mathcal{P}_{2,5,0}$ has $\binom{5}{2} + \binom{5}{2}= 20$
elements. It has five incomparable pairs
$\langle ij \rangle \langle kl \rangle $,
five incomparable pairs $[ ij ] [ kl ] $,
and $25$ mixed incomparable pairs
$\langle ij \rangle [ kl]$.
These are the initial monomials of the $35$ ideal generators
 in Example \ref{ex:twofivezero}.
\end{example}

The following result concerns the number of incomparable pairs 
in the poset $\Pcal_{k,n,r}$.

\begin{lemma} \label{lem:mixed}
    The number of mixed incomparable pairs $\langle i_1, \dots, i_k \rangle [j_1, \dots, j_k]$
      is $\binom{n}{k-r-1}^2$.
\end{lemma}

\begin{proof}
Let $\lambda = (n-k+r, \ k)$ and $\mu = (r)$. The generating function for semi-standard skew tableaux of shape $\lambda / \mu$ is the {\em skew Schur polynomial} $s_{\lambda / \mu}$, which can be written as 
\[
    s_{\lambda / \mu} \,\,=\,\, \sum_{\nu} c_{\mu, \nu}^{\lambda} s_{\nu}.
\]
Here $c^{\lambda}_{\mu, \nu}$ are the  Littlewood-Richardson coefficients. 
In our special case $\mu = (r)$, Pieri's rule tells us that
   $c^{\lambda}_{\mu, \nu} = 1$ if $\nu$ is obtained from $\lambda$ by removing $r$ boxes and $0$ otherwise.           Using the hook-content formula to evaluate each individual $s_\nu(1,\dots, 1)$, we
   can now compute the number of semi-standard skew tableaux of shape $\lambda / \mu$ with fillings in $[n]$.
   That number is
\begin{equation} \label{eq:SSSYT}  \begin{matrix}
    s_{\lambda/\mu}(1,1,\ldots,1) \quad \,
    & = & \,\,  \sum_{ \substack{\nu \text{ obtained from $\lambda$ } \\ \text{by removing $r$ boxes} } } s_{\nu}(1,1, \dots,1) \medskip \\ 
    &=& \quad \sum_{\ell = 0}^{r} s_{(n-k+\ell, k - \ell)}(1,1, \dots, 1) \smallskip \\
    &=&\quad \sum_{\ell = 0}^{r}  \left[ \binom{n}{k-\ell}^2  - \binom{n}{k-\ell-1}^2 \right ]\smallskip \\
    &=&\quad \binom{n}{k}^{\! 2} -\, \binom{n}{k-r-1}^{\! 2} . 
\end{matrix}
\end{equation}
Since  $\binom{n}{k}^2$ is the number of all tableaux,
we see that the number of non-standard skew tableau of shape 
$\lambda/\mu$ equals $\binom{n}{k-r-1}^2$.
The equivalence in (\ref{eq:3equivalent}) now completes the proof.
\end{proof}

Theorem \ref{thm:main1} and Lemma \ref{lem:mixed}
imply the following result about the spinor-helicity variety.

\begin{corollary} \label{cor:arequadratic}
The number of minimal generators of $I_{k,n,r}$
equals twice (\ref{eq:hookcontent}) plus $\binom{n}{k-r-1}^2$.
These generators are quadratic and we can 
arrange them to form a
reduced Gr\"obner basis.
\end{corollary}

After its dimension, the second-most important invariant of a variety in a projective space is its
degree. For a variety in a product of two projective spaces, one considers
the bidegree, which is a homogeneous polynomial in two variables $s$ and $t$.
The degree of that polynomial is the codimension of the variety.
We saw an example in (\ref{eq:cohom250}). By definition,
the {\em bidegree} of ${\rm SH}(k,n,r)$ is its class in the cohomology ring
of ambient product of projective spaces:
$$ H^*\left(\PP^{\binom{n}{k}-1} \times \PP^{\binom{n}{k}-1}, \ZZ \right) \,\, = \,\,
\ZZ[s,t]\,/\bigl\langle\, s^{\binom{n}{k}},\, t^{\binom{n}{k}}\, \bigr\rangle. $$

We now present a general formula for the cohomology class of
the spinor-helicity variety.

\begin{corollary} \label{cor:bidegree}
The bidegree of $\,{\rm SH}(k,n,r)$ is equal to
\begin{equation}
\label{eq:bidegree} (st)^{\binom{n}{k}-k(n-k)-1} \cdot \sum 
c(i_1 i_2 \cdots i_k) \cdot c(j_1 j_2 \cdots j_k) \cdot
s^{i_1+i_2+\cdots+i_k - \binom{k+1}{2}} \cdot t^{j_1+j_2+\cdots+j_k - \binom{k+1}{2}},
\end{equation}
where we sum  over all covering relations in (\ref{eq:coveringrel}),
and $c(i_1 i_2 \cdots i_k)$ denotes the number of maximal chains
from $\langle i_1 i_2 \cdots i_k \rangle$ to the top element
$\langle n{-}k{+}1 \,\cdots \, n{-}1 \, n \rangle$ in
Young's lattice~$Y_{k,n}$.
\end{corollary}

\begin{proof}
The bidegree of $I_{k,n,r}$ equals the bidegree of the
initial monomial ideal ${\rm in}(I_{k,n,r})$.
The latter is generated by the
incomparable pairs in $\mathcal{P}_{k,n,r}$. The bidegree is
the multidegree of \cite[\S 8.5]{CCA} for the
$\ZZ^2$-grading at hand. It is additive over  top-dimensional
 primary components, and we can use the formula
in \cite[Theorem 8.44]{CCA} for its evaluation.
The associated primes of ${\rm in}(I_{k,n,r})$ correspond to the
maximal chains of $\mathcal{P}_{k,n,r}$. Each maximal chain
starts out at the bottom of $\tilde Y_{k,n}$, it
uses precisely one of the covering relations 
in (\ref{eq:coveringrel}) to transition from $\tilde Y_{k,n}$
to $Y_{k,n}$, and it then proceeds to the top of $Y_{k,n}$.
Thus there are precisely 
$c(i_1 i_2 \cdots i_k) \cdot c(j_1 j_2 \cdots j_k) $
maximal chains which use the
specific covering relation in (\ref{eq:coveringrel}).
The associated monomial in $s$ and $t$
records the height at which the transition 
from $\tilde Y_{k,n}$ to $Y_{k,n}$ in $\mathcal{P}_{k,n,r}$
takes place.
\end{proof}

\begin{example}[$k=2,n=6,r=0$]
The bidegree of $I_{2,6,0}$ is a sum of
monomials $s^i t^j$ over the maximal chains of the poset
$\mathcal{P}_{2,6,0}$
 in Figure \ref{fig:P_(2,6,0)}.
 The degree $i+j=16$ of each monomial is the
 codimension of ${\rm SH}(2,6,0)$
 in $\PP^{14} \times \PP^{14}$.
 Counting paths in Young's lattice $Y_{2,6}$,   we see
$$ c(12) = 14,\,
c(13) = 14,\,
c(14) = 9,\,
c(23) =  5,\,
c(24) = 5,\,
c(34) = 2.
$$
The bidegree of $I_{2,6,0}$ is the polynomial
$\,28 s^6 t^{10} +
70 s^7 t^9 +
90 s^8 s^8 + 
70 s^9 t^7
+ 28 s ^{10} t^6$
$$ 
=\, c(12) c(34) s^6 t^{10} +
c(13) c(24) s^7 t^9 +
2 c(14) c(23) s^8 s^8
+ c(24) c(13) s^9 t^7
+ c(34) c(12) s ^{10} t^6 .
$$
Hence total number of maximal chains in $\mathcal{P}_{2,6,0}$
equals $28+70+90+70+28 = 286$.
\end{example}

\begin{example}[$k=3,n=7,r=1$]
The poset  $\mathcal{P}_{3,7,1}$ has $70$ elements and
$312816$ maximal chains.
It arises from $Y_{3,7}$ and $\tilde Y_{3,7}$
by adding six covering relations, as shown in Figure \ref{fig:P_(3,7,1)}.
The ideal $I_{3,7,1}$ has $140+140+49 = 329$
minimal generators, one for each incomparable pair;
by Corollary \ref{cor:arequadratic}.
The blue numbers $c(ijk)$ count 
maximal chains from $\langle ijk \rangle $ to $\langle 567 \rangle$
in $Y_{3,7}$ or maximal chains
from $[567]$ to $[ijk]$ in $\tilde Y_{3,7}$.
By Corollary \ref{cor:bidegree},
the bidegree  of $I_{3,7,1}$ equals    
    $(st)^{22} \bigl( 462 \cdot 56 \,s^4 
    \,+\, 462 \cdot 168 \,s^3t \,+\, (252 \cdot 210 + 210 \cdot 252) \,s^2t^2 
    \,+\, 168 \cdot 462\, st^3     \,+\,  56 \cdot 462 \,t^4 \bigr)$.
\end{example}

\begin{figure}[!htbp]
    \centering
    \begin{center}
\scalebox{0.6}{
    \begin{tikzpicture}
            \draw (4,8) node {$ \langle 167 \rangle $};
            \draw[blue] (3.5,8.20) node {\small 1};
            \draw (2,6) node { $ \langle 157 \rangle $};
            \draw[blue] (1.5,6.20) node {\small 4};
            \draw (0,4) node {$ \langle  147 \rangle $};
            \draw[blue] (-0.5,4.20) node {\small 9};
            \draw (-2,2) node { $ \langle 137 \rangle $};
            \draw[blue] (-2.7,2.20) node {\small 14};
            \draw (-4,0) node {$ \langle 127 \rangle $};
            \draw[blue] (-4.7,0.20) node {\small 14};

            \draw (-3.75,0.25)--(-2.25,1.75);
            \draw (-1.75,2.25)--(-0.25,3.75);
            \draw (0.25,4.25)--(1.75,5.75);
            \draw (2.25,6.25)--(3.75,7.75);

            \draw (-3.75,-0.25)--(-2.25,-1.75);
            \draw (-1.75,1.75)--(-0.25,0.25);
            \draw (0.25,3.75)--(1.75,2.25);
            \draw (2.25,5.75)--(3.75,4.25);

            \draw (4,4) node {$ \langle 156 \rangle $};
            \draw[blue] (3.5,4.20) node {\small 10};
            \draw (2,2) node {$ \langle 146 \rangle $};
            \draw[blue] (1.4,2.20) node {\small 35};
            \draw (0,0) node {$ \langle 136 \rangle $};
            \draw[blue] (-0.7,0.20) node {\small 70};
            \draw (-2,-2) node {$ \langle 126 \rangle $};
            \draw[blue] (-2.7,-1.80) node {\small 84};

            \draw (-1.75,-1.75)--(-0.25,-0.25);
            \draw (0.25,0.25)--(1.75,1.75);
            \draw (2.25,2.25)--(3.75,3.75);

            \draw (-1.75,-2.25)--(-0.25,-3.75);
            \draw (0.25,-0.25)--(1.75,-1.75);
            \draw (2.25,1.75)--(3.75,0.25);

            \draw (4,0) node {$ \langle 145 \rangle $};
            \draw[blue] (3.4,0.20) node {\small 56};
            \draw (2,-2) node {$ \langle 135 \rangle $};
            \draw[blue] (1.2,-1.80) node {\small 168};
            \draw (0,-4) node {$ \langle 125 \rangle $};
            \draw[blue] (-0.8,-3.80) node {\small 252};

            \draw (0.25,-3.75)--(1.75,-2.25);
            \draw (2.25,-1.75)--(3.75,-0.25);

            \draw (0.25,-4.25)--(1.75,-5.75);
            \draw (2.25,-2.25)--(3.75,-3.75);

            \draw (4,-4) node {$ \langle 134 \rangle $};
            \draw[blue] (3.2,-3.80) node {\small 210};
            \draw (2,-6) node {$ \langle 124 \rangle $};
            \draw[blue] (1.2,-5.80) node {\small 462};

            \draw (2.25,-5.75)--(3.75,-4.25);
            
            \draw (2.25,-6.25)--(3.75,-7.75);
        
            \draw (4,-8) node {$ \langle 123 \rangle $};
            \draw[blue] (3.2,-7.80) node {\small 462};

            {\color{red}

            
            \draw[black] (4,9) node {$ \langle 267 \rangle $};
            \draw[black] (4,8.2)--(4,8.8);
            \draw[blue] (3.5,9.20) node {\small 1};
            
            \draw[black] (2,7) node {$ \langle 257 \rangle $};
            \draw[black] (2,6.2)--(2,6.8);
            \draw[blue] (1.5,7.20) node {\small 3};
            
            \draw[black] (0,5) node {$ \langle 247 \rangle $};
            \draw[black] (0,4.2)--(0,4.8);
            \draw[blue] (-0.5,5.20) node {\small 5};
            
            \draw[black] (-2,3) node {$ \langle 237 \rangle $};
            \draw[black] (-2,2.2)--(-2,2.8);
            \draw[blue] (-2.5,3.20) node {\small 5};

            \draw (-1.75,3.25)--(-0.25,4.75);
            \draw (0.25,5.25)--(1.75,6.75);
            \draw (2.25,7.25)--(3.75,8.75);

            \draw (-1.75,2.75)--(-0.25,1.25);
            \draw (0.25,4.75)--(1.75,3.25);
            \draw (2.25,6.75)--(3.75,5.25);

            \draw[black] (4,5) node {$ \langle 256 \rangle $};
            \draw[black] (4,4.2)--(4,4.8);
            \draw[blue] (3.5,5.20) node {\small 6};
            
            \draw[black] (2,3) node {$ \langle 246 \rangle $};
            \draw[black] (2,2.2)--(2,2.8);
            \draw[blue] (1.4,3.20) node {\small 16};
            
            \draw[black] (0,1) node {$ \langle 236 \rangle $};
            \draw[black] (0,0.2)--(0,0.8);
            \draw[blue] (-0.7,1.20) node {\small 21};

            \draw (0.25,1.25)--(1.75,2.75);
            \draw (2.25,3.25)--(3.75,4.75);

            \draw (0.25,0.75)--(1.75,-0.75);
            \draw (2.25,2.75)--(3.75,1.25);

            \draw[black] (4,1) node {$ \langle 245 \rangle $};
            \draw[black] (4,0.2)--(4,0.8);
            \draw[blue] (3.4,1.20) node {\small 21};
            
            \draw[black] (2,-1) node {$ \langle 235 \rangle $};
            \draw[black] (2,-1.8)--(2,-1.2);
            \draw[blue] (1.3,-0.8) node {\small 42};

            \draw (2.25,-0.75)--(3.75,0.75);
            
            \draw (2.25,-1.25)--(3.75,-2.75);
        
            \draw[black] (4,-3) node {$ \langle 234 \rangle $};
            \draw[black] (4,-3.8)--(4,-3.2);
            \draw[blue] (3.3, -2.80) node {\small 42};
            
            }

         
            {\color{blue}
            
            \draw[black] (4,10) node {$ \langle 367 \rangle $};
            \draw[black] (4,9.2)--(4,9.8);
            \draw[blue] (3.5,10.20) node {\small 1};
            
            \draw[black] (2,8) node {$ \langle 357 \rangle $};
            \draw[black] (2,7.2)--(2,7.8);
            \draw[blue] (1.5,8.20) node {\small 2};
            
            \draw[black] (0,6) node {$ \langle 347 \rangle $};
            \draw[black] (0,5.2)--(0,5.8);
            \draw[blue] (-0.5,6.20) node {\small 2};

            \draw (0.25,6.25)--(1.75,7.75);
            \draw (2.25,8.25)--(3.75,9.75);

            \draw (0.25,5.75)--(1.75,4.25);
            \draw (2.25,7.75)--(3.75,6.25);
    
            \draw[black] (4,6) node {$ \langle 356 \rangle $};
            \draw[black] (4,5.2)--(4,5.8);
            \draw[blue] (3.5,6.20) node {\small 3};
            
            \draw[black] (2,4) node {$ \langle 346 \rangle $};
            \draw[black] (2,3.2)--(2,3.8);
            \draw[blue] (1.5,4.20) node {\small 5};

            \draw (2.25,4.25)--(3.75,5.75);

            \draw (2.25,3.75)--(3.75,2.25);
            
            \draw[black] (4,2) node {$ \langle 345 \rangle $};
            \draw[black] (4,1.2)--(4,1.8);
            \draw[blue] (3.5,2.20) node {\small 5};
            }

            {\color{green}
            
            \draw[black] (4,11) node {$ \langle 467 \rangle $};
            \draw[black] (4,10.2)--(4,10.8);
            \draw[blue] (3.5,11.20) node {\small 1};
            
            \draw[black] (2,9) node {$ \langle 457 \rangle $};
            \draw[black] (2,8.2)--(2,8.8);
            \draw[blue] (1.5,9.20) node {\small 1};
               
            \draw (2.25,9.25)--(3.75,10.75);
      
            \draw (2.25,8.75)--(3.75,7.25);
    
            \draw[black] (4,7) node {$ \langle 456 \rangle $};
            \draw[black] (4,6.2)--(4,6.8);
            \draw[blue] (3.5,7.20) node {\small 1};
            }

            \draw[black] (4,12) node {$ \langle 567 \rangle $};
            \draw[black] (4,11.2)--(4,11.8);
            \draw[blue] (3.5,12.2) node {\small 1};

        \begin{scope}[shift={(15,-8)} , rotate=180]
            \draw (4,8) node {$[167]$};
            \draw[blue] (3.5,7.80) node {\small 1};
            \draw (2,6) node { $[157]$};
            \draw[blue] (1.5,5.80) node {\small 4};
            \draw (0,4) node {$[ 147]$};
            \draw[blue] (-0.5,3.80) node {\small 9};
            \draw (-2,2) node { $[137]$};
            \draw[blue] (-2.6,1.80) node {\small 14};
            \draw (-4,0) node {$[127]$};
            \draw[blue] (-4.6,-0.20) node {\small 14};

            \draw (-3.75,0.25)--(-2.25,1.75);
            \draw (-1.75,2.25)--(-0.25,3.75);
            \draw (0.25,4.25)--(1.75,5.75);
            \draw (2.25,6.25)--(3.75,7.75);

            \draw (-3.75,-0.25)--(-2.25,-1.75);
            \draw (-1.75,1.75)--(-0.25,0.25);
            \draw (0.25,3.75)--(1.75,2.25);
            \draw (2.25,5.75)--(3.75,4.25);

            \draw (4,4) node {$[156]$};
            \draw[blue] (3.4,3.80) node {\small 10};
            \draw (2,2) node {$[146]$};
            \draw[blue] (1.3,1.80) node {\small 35};
            \draw (0,0) node {$[136]$};
            \draw[blue] (-0.6,-0.20) node {\small 70};
            \draw (-2,-2) node {$[126]$};
            \draw[blue] (-2.6,-2.2) node {\small 84};

            \draw (-1.75,-1.75)--(-0.25,-0.25);
            \draw (0.25,0.25)--(1.75,1.75);
            \draw (2.25,2.25)--(3.75,3.75);

            \draw (-1.75,-2.25)--(-0.25,-3.75);
            \draw (0.25,-0.25)--(1.75,-1.75);
            \draw (2.25,1.75)--(3.75,0.25);

            \draw (4,0) node {$[145]$};
            \draw[blue] (3.4,-0.20) node {\small 56};
            \draw (2,-2) node {$[135]$};
            \draw[blue] (1.3,-2.2) node {\small 168};
            \draw (0,-4) node {$[125]$};
            \draw[blue] (-0.7,-4.20) node {\small 252};

            \draw (0.25,-3.75)--(1.75,-2.25);
            \draw (2.25,-1.75)--(3.75,-0.25);

            \draw (0.25,-4.25)--(1.75,-5.75);
            \draw (2.25,-2.25)--(3.75,-3.75);

            \draw (4,-4) node {$[134]$};
            \draw[blue] (3.3,-4.20) node {\small 210};
            \draw (2,-6) node {$[124]$};
            \draw[blue] (1.3,-6.20) node {\small 462};

            \draw (2.25,-5.75)--(3.75,-4.25);
            
            \draw (2.25,-6.25)--(3.75,-7.75);
        
            \draw (4,-8) node {$[123]$};
            \draw[blue] (3.3,-8.20) node {\small 462};

            {\color{red}

            
            \draw[black] (4,9) node {$[267]$};
            \draw[black] (4,8.2)--(4,8.8);
            \draw[blue] (3.5,8.80) node {\small 1};
            
            \draw[black] (2,7) node {$[257]$};
            \draw[black] (2,6.2)--(2,6.8);
            \draw[blue] (1.5,6.80) node {\small 3};
            
            \draw[black] (0,5) node {$[247]$};
            \draw[black] (0,4.2)--(0,4.8);
            \draw[blue] (-0.5,4.80) node {\small 5};
            
            \draw[black] (-2,3) node {$[237]$};
            \draw[black] (-2,2.2)--(-2,2.8);
            \draw[blue] (-2.5,2.80) node {\small 5};

            \draw (-1.75,3.25)--(-0.25,4.75);
            \draw (0.25,5.25)--(1.75,6.75);
            \draw (2.25,7.25)--(3.75,8.75);

            \draw (-1.75,2.75)--(-0.25,1.25);
            \draw (0.25,4.75)--(1.75,3.25);
            \draw (2.25,6.75)--(3.75,5.25);

            \draw[black] (4,5) node {$[256]$};
            \draw[black] (4,4.2)--(4,4.8);
            \draw[blue] (3.5,4.80) node {\small 6};
            
            \draw[black] (2,3) node {$[246]$};
            \draw[black] (2,2.2)--(2,2.8);
            \draw[blue] (1.4,2.80) node {\small 16};
            
            \draw[black] (0,1) node {$[236]$};
            \draw[black] (0,0.2)--(0,0.8);
            \draw[blue] (-0.6,0.80) node {\small 21};

            \draw (0.25,1.25)--(1.75,2.75);
            \draw (2.25,3.25)--(3.75,4.75);

            \draw (0.25,0.75)--(1.75,-0.75);
            \draw (2.25,2.75)--(3.75,1.25);

            \draw[black] (4,1) node {$[245]$};
            \draw[black] (4,0.2)--(4,0.8);
            \draw[blue] (3.4,0.80) node {\small 21};
            
            \draw[black] (2,-1) node {$[235]$};
            \draw[black] (2,-1.8)--(2,-1.2);
            \draw[blue] (1.4,-1.2) node {\small 42};

            \draw (2.25,-0.75)--(3.75,0.75);
            
            \draw (2.25,-1.25)--(3.75,-2.75);
        
            \draw[black] (4,-3) node {$[234]$};
            \draw[black] (4,-3.8)--(4,-3.2);
            \draw[blue] (3.4,-3.20) node {\small 42};
            
            }

         
            {\color{blue}
            
            \draw[black] (4,10) node {$[367]$};
            \draw[black] (4,9.2)--(4,9.8);
            \draw[blue] (3.5,9.80) node {\small 1};
            
            \draw[black] (2,8) node {$[357]$};
            \draw[black] (2,7.2)--(2,7.8);
            \draw[blue] (1.5,7.80) node {\small 2};
            
            \draw[black] (0,6) node {$[347]$};
            \draw[black] (0,5.2)--(0,5.8);
            \draw[blue] (-0.5,5.80) node {\small 2};

            \draw (0.25,6.25)--(1.75,7.75);
            \draw (2.25,8.25)--(3.75,9.75);

            \draw (0.25,5.75)--(1.75,4.25);
            \draw (2.25,7.75)--(3.75,6.25);
    
            \draw[black] (4,6) node {$[356]$};
            \draw[black] (4,5.2)--(4,5.8);
            \draw[blue] (3.5,5.80) node {\small 3};
            
            \draw[black] (2,4) node {$[346]$};
            \draw[black] (2,3.2)--(2,3.8);
            \draw[blue] (1.5,3.80) node {\small 5};

            \draw (2.25,4.25)--(3.75,5.75);

            \draw (2.25,3.75)--(3.75,2.25);
            
            \draw[black] (4,2) node {$[345]$};
            \draw[black] (4,1.2)--(4,1.8);
            \draw[blue] (3.5,1.80) node {\small 5};
            }

            {\color{green}
            
            \draw[black] (4,11) node {$[467]$};
            \draw[black] (4,10.2)--(4,10.8);
            \draw[blue] (3.5,10.80) node {\small 1};
            
            \draw[black] (2,9) node {$[457]$};
            \draw[black] (2,8.2)--(2,8.8);
            \draw[blue] (1.5,8.80) node {\small 1};
               
            \draw (2.25,9.25)--(3.75,10.75);
      
            \draw (2.25,8.75)--(3.75,7.25);
    
            \draw[black] (4,7) node {$[456]$};
            \draw[black] (4,6.2)--(4,6.8);
            \draw[blue] (3.5,6.80) node {\small 1};
            }

            \draw[black] (4,12) node {$[567]$};
            \draw[black] (4,11.2)--(4,11.8);
            \draw[blue] (3.5,11.80) node {\small 1};
        \end{scope}


        \draw[black, dashed] (0.25,-4) .. controls (4,-5) and (6,-5) .. (10.75,-4);
        \draw[black, dashed] (4.25,-4) .. controls (8,-3) and (10,-3) .. (15,-4);
        \draw[black, dashed] (2.25,-2)--(12.75,-2);
        \draw[black, dashed] (4.25,0)--(10.75,0);
        \draw[black, dashed] (2.25,-6)--(12.75,-6);
        \draw[black, dashed] (4.25,-8)--(10.75,-8);

    \end{tikzpicture}
    }
\end{center}
    \caption{ The poset $\Pcal_{3,7,1}$ governs the 
combinatorics of the
        variety ${\rm SH}(3,7,1) \subset \PP^{34} \times \PP^{34}$.
    \label{fig:P_(3,7,1)} }
\end{figure}

\section{Bilinear Relations and Khovanskii Bases}
\label{sec3}

In Section \ref{sec2} we took a route into
the combinatorial commutative algebra of the 
spinor-helicity variety. This journey continues
in this section. We begin by taking a closer look
at the bilinear equations in the ideal $I_{k,n,r}$.
Thereafter, we introduce a toric degeneration of the variety
${\rm SH}(k,n,r)$, based on the poset  $\mathcal{P}_{k,n,r}$.
The associated Khovanskii basis 
of the coordinate ring matches our earlier Gr\"obner basis.

We now introduce two matrices $P$ and $Q$
whose rows are indexed by $\binom{[n]}{k-r-1}$
and whose columns are indexed by $\binom{[n]}{r+1}$.
Here $[n] = \{1,2,\ldots,n\}$ and 
$\binom{[n]}{s}$ is the set of subsets of size $s$ in $[n]$.
The entries of our matrices are given by concatenating row labels
and column labels
\[
    P_{I,J} \,=\,  \langle \,I \, J \,\rangle  
\quad \text{and }  \quad  Q_{I,J} \,=\,  [\,I\,J\,]
\qquad \hbox{for $\,I \in  \binom{[n]}{k-r-1}\,$
and $\,J \in \binom{[n]}{r+1}$.}
\]
Here $I$ and $J$ are increasing sequences
which we concatenate.
We pass to the sorted
Pl\"ucker coordinates $\,\langle I \,\cup \, J  \rangle$
and $\,[I \, \cup \, J]$ from Section \ref{sec2}
by multiplying with $-1,+1$ or $0$, as  in (\ref{eq:sorttt}).
In particular, this means that
$\,\langle I \, J  \rangle \, = \,0\,$
and $\,[I \, J] \,  = \,0\,$
whenever $I \cap J \not= \emptyset$.

\begin{example}[$k=4, n=5, r=1$] The matrix $P$ is square of format $10 \times 10$.
The rows and columns of $P$ are labeled by
$ 12, 13, 14, 15, 23, 24, 25, 34, 35,45$, in this order. We find
\[
\resizebox{\textwidth}{!}{$P\, =\,  \begin{pmatrix}
           0  &  0  &  0  &  0  &  0 &  0 & 0 & \langle1234\rangle  & \langle1235\rangle & \langle1245\rangle\\
            0  &  0  &  0  &  0  &  0 &  \langle1324\rangle & \langle1345\rangle &  0  & 0 & \langle1345\rangle\\
            0  &  0  &  0  &  0  &  \langle1423\rangle &  0 &  \langle1425\rangle &  0  & \langle1435\rangle & 0\\
            0  &  0  &  0  &  0  &  \langle1523\rangle & \langle1524\rangle & 0 &  \langle1534\rangle  & 0 & 0\\
           0  &  0  &  \langle2314\rangle  &  \langle2315\rangle  &  0 &  0 & 0 &  0  & 0 & \langle2345\rangle\\
            0  &  \langle2413\rangle  &  0  &  \langle2415\rangle  &  0 &  0 & 0 & 0 & \langle2435\rangle & 0\\
            0  &  \langle2513\rangle  &  \langle2514\rangle  &  0  &  0 &  0 & 0 &  \langle2534\rangle  & 0 & 0\\
            \langle3412\rangle  &  0  &  0  &  \langle3415\rangle  &  0 &  0 & \langle3425\rangle &  0  & 0 & 0\\
            \langle3512\rangle  &  0  &  \langle3514\rangle  &  0 &  0 &  \langle3524\rangle & 0 &  0  & 0 & 0\\
            \langle4512\rangle  &  \langle4513\rangle  &  0  &  0  &  \langle4523\rangle &  0 & 0 &  0  & 0 & 0\\
          \end{pmatrix}\!. 
          $}
\]
Each entry is now replaced by Pl\"ucker coordinates with increasing indices.
For instance, we replace  $\langle 1324 \rangle $ by
$- \langle    1234 \rangle $. The matrix $Q$ is identical to $P$ but
with square brackets $[ijkl]$.
\end{example}

\begin{example}[$k=2,r=0$]
This is the case of most interest in physics.
Here $P$ and $Q$ are the skew-symmetric 
$n \times n$ matrices of Pl\"ucker coordinates,
 shown for $n=5$ in Example~\ref{ex:twofivezero}.
\end{example}

\begin{example}[$r=k-1$]
The variety ${\rm SH}(k,n,k-1)$ is
a hypersurface in ${\rm Gr}(k,n) \times {\rm Gr}(k,n)$.
Here $P$ and $Q$ are the row vectors of length $\binom{n}{k}$
whose entries are the  Pl\"ucker coordinates. 
The defining equation of this hypersurface is the inner product of the two  Pl\"ucker vectors:
$$ P Q^T \, = \,\sum \langle i_1 i_2 \cdots i_k \rangle [i_1 i_2 \cdots i_k ] \,\, = \,\, 0 . $$
\end{example}

These examples guide us towards our next theorem, which is the main result in Section~\ref{sec3}.

\begin{theorem} \label{thm:PQ}
The entries of the matrix $\,P Q^T$ generate the prime ideal of
the spinor-helicity variety ${\rm SH}(k,n,r)$ in the coordinate ring of $\,{\rm Gr}(k,n) \times {\rm Gr}(k,n)$.
In symbols, we have
\begin{equation} \label{eq:idealeq}
    I_{k,n,r} \,\,= \,\,J_{n,k} + \tilde J_{k,n} + \langle \text{entries of } P Q^T \rangle.
\end{equation}
\end{theorem}

\begin{proof}
We first show that the entries of $P Q^T$ vanish on the spinor-helicity variety.
Fix any point in ${\rm SH}(k,n,r)$, represented by a pair
of $k \times n$ matrices $\lambda$ and $\tilde \lambda$ such that
$\lambda \cdot \tilde \lambda^T$ has rank $\leq r$.
By passing to the $(r+1)$-st exterior power, we find that
$\, \wedge_{r+1} \lambda \cdot (\wedge_{r+1} \tilde \lambda)^T\,$
is the zero matrix of format $\binom{k}{r+1} \times \binom{k}{r+1}$.
In other words, the row spaces of the $\binom{k}{r+1} \times \binom{n}{r+1}$
matrices $\wedge_{r+1} \lambda$ and $\wedge_{r+1} \tilde \lambda$
are orthogonal to each other. The row vectors of the matrices $P$ and $Q$
are elements in these row spaces. Therefore $P \cdot Q^T = 0$
holds for our point $(\lambda,\tilde \lambda)$.

The previous paragraph shows that the right hand side of  (\ref{eq:idealeq})
is contained in the left hand side. Both ideals are generated by quadrics,
and they contain the Pl\"ucker ideals $J_{k,n}$ and $\tilde J_{k,n}$.
It therefore suffices to show that the entries of $PQ^T$ span
the space of all bilinear quadrics in the ideal $I_{k,n,r}$.
We know from Lemma \ref{lem:mixed} and
Corollary \ref{cor:arequadratic} that this
space has dimension $\binom{n}{k-r-1}^2$.
This number coincides with the number of entries in the square matrix $PQ^T$.
It therefore suffices to show that the entries of $PQ^T$ are linearly independent
over $\CC$.

We shall prove this by contradiction.
The entry of $PQ^T$ in row $I$ and column $J$ equals
 $$ f_{IJ} \,\,\,= \sum_{L \in \binom{[n]}{r+1}} \epsilon_{I, L}\,\epsilon_{J,L} \ \langle IL \rangle  [JL],
 $$
 where $\epsilon_{I,L}=\pm 1$ is the sign of the permutation that sorts 
 the string $IL$. Suppose that$$
    \sum_{I,J \in \binom{[n]}{k-r-1}} \!\!\!\! \alpha_{I,J} \cdot f_{IJ} \,\,= \,\,0 \qquad
\hbox{for some scalars $\alpha_{IJ} \in \CC$.}
$$
We must show that each $\alpha_{IJ}$ is zero.
        The previous equation can be rewritten as follows:
 $$       \sum_{ \substack{I', J' \in \binom{[n]}{k} \\ |I' \cap J'| \geq r+1 }}      \biggl[\,\,
        \sum_{\substack{ L \subset I' \cap J' \\ |L|=r+1} }\!\! \epsilon_{I'\setminus L , L}  \,\epsilon_{J'\setminus L, L} \ \alpha_{I' \setminus L, J' \setminus L} 
       \, \biggr]\, \langle \,I'\, \rangle \,[\,J'\, ] \,\,\,= \,\,\,0. $$
    From this we conclude that,  for any two $k$-subsets $I',J'$ with $|I' \cap J'|\geq r+1$, we have
\begin{equation}
\label{eq:weconclude}    \begin{matrix} \qquad
     \sum_{\substack{L \in \binom{ I' \cap J'}{r+1}}}
         \epsilon_{I'\setminus L , L} \,\epsilon_{J'\setminus L, L} \ \alpha_{I' \setminus L, J' \setminus L} \,\,=\,\, 0
        \qquad \hbox{for all} \,\,\, I',J' \in \binom{ [n]}{k}.  \end{matrix}
\end{equation}        

    Our goal is to show that all $\alpha$'s are zero.
    First consider the case $I' = J'$. Here (\ref{eq:weconclude})    reads
    \[        \sum_{L \subset I',  |L| = r+1} \!\!\!\! \alpha_{I' \setminus L,  I' \setminus L} 
        \,\,= \,\,0.    \] 
   We  write these equations in the form $B a = 0$ where $a_I = \alpha_{I,I}$ and   $B$ is an
        $\binom{n}{k} \times \binom{n}{k-r-1}$ matrix  with entries in $\{0,1\}$.     
        The row indices are subsets $ J \in \binom{[n]}{k} $ and the column indices are
        subsets $  I \subset \binom{[n]}{k-r-1}$.
   The matrix entry $B_{J, I} $ equals $1$ when $I \subset J$
        and it is $0$ otherwise. 

It is a known result in combinatorics that the columns of the matrix $B$
are linearly independent. The context is the spectral theory of the
{\em Johnson graph}, which is developed in \cite[Chapter 6]{GM16}.
To be precise, the desired identity
$\,{\rm rank}(B) = \binom{n}{k-r-1}\,$ can be found in
 \cite[Theorem 6.3.3]{GM16}.
           From this statement we deduce that
            $\alpha_{I,I} = 0$ for all $I \in \binom{[n]}{k-r-1}$. 

      For the general case, we fix
      $I' = I_{0} \sqcup K$ and $J' = J_0 \sqcup K$
      where $I_0, J_0 \in \binom{[n]}{s}$ are disjoint,
      with $s \geq 1$,
      and  $K \subset [n] \backslash (I_0 \sqcup J_0)$ has size $k-s $.
         From  \eqref{eq:weconclude} we obtain the equations
      \begin{equation}\label{eq:geneLinIndep} 
      \!\!\! \sum_{\substack{L \subset K \\ |L|=r+1}}\!\! \epsilon_{(I_0 \sqcup K) \setminus L,L} \cdot \epsilon_{(J_0 \sqcup K) \setminus L, L} \cdot
          \alpha_{(I_0 \sqcup K) \setminus L,   (J_0 \sqcup K) \setminus L} 
          \,=\, 0 \quad \text{for any } K \in \binom{[n]\backslash (I_0 \sqcup J_0)}{k-s}.
      \end{equation}
      Suppressing $I_0$ and $J_0$ from the indices of $\alpha$, we rewrite
      (\ref{eq:geneLinIndep})
        in matrix form $U B V a = 0$,~where
      \begin{itemize}
          \item $a_{I} = \alpha_{I_0 \sqcup I, J_0 \sqcup I}$ for any 
          subset $I \in \binom{[n] \setminus (I_0 \sqcup J_0)}{k-s-r-1}$. 
          To get to \eqref{eq:geneLinIndep}, we would set $I = K \backslash L$.
      \item $B$ is a $\binom{n-2s}{k - s} \times \binom{n-2s}{k-s-r-1}$ matrix with entries in $\{0,1\}$. The columns of $B$ are indexed by subsets $I$ of size $k{-}s{-}r{-}1$ of $[n] \backslash 
      (I_0 \sqcup J_0)$ and the rows are indexed by
      subsets $K$ of size $k{-}s$ subsets of $[n] \backslash (I_0 \sqcup J_0)$.
      The entries are  $B_{K,I} = 1$ of $ I \subset K$ and $0$ otherwise.
          \item $U$ and $V$ are diagonal matrices  of size $\binom{n-2s}{k-s} $ 
          and $\binom{n-2s}{k-s-r-1} $ respectively, with entries
          \[
            U_{K,K} \,=\, (-1)^{N}
            \quad \text{ and } \quad V_{I,I} \,=\, (-1)^{M},
          \]
          where $N$ counts the elements in $I_0 \sqcup J_0$ that are larger than elements in $K$ and $M$ counts the elements in $I_0 \sqcup J_0$ that are larger than elements in $I$.
          In symbols,
          \[
          N \,\,= \sum_{i \in I_0 \sqcup J_0} \sum_{\ell \in K} 1_{i > \ell} \qquad \text{and} \qquad 
          M\,\, = \sum_{i \in I_0 \sqcup J_0} \sum_{\ell \in I} 1_{i > \ell}.
          \vspace{-0.2cm} \]
      \end{itemize}      
       Writing $\,I \subset K\,$ and $\,L = K \backslash I$, we find $\,\epsilon_{(I_0 \sqcup K) \setminus L,L} \cdot \epsilon_{(J_0 \sqcup K) \setminus L, L}  = (-1)^{N+M}\,$ for the sign in 
       
      \medskip
\noindent \eqref{eq:geneLinIndep}.
      Again, by virtue of \cite[Theorem 6.3.3]{GM16}, we have ${\rm rank}(B) = \binom{n-2s}{k-s-r-1}$ and hence $a_{I} = \alpha_{I_0 \sqcup I, J_0 \sqcup I} = 0$ for any set 
      $I \subset [n] \backslash (I_0 \sqcup J_0)$ of size $k{-}s{-}r{-}1$. 
      Since this holds  for any pair of disjoint index sets $I_0,J_0
      \in \binom{[n]}{s}$ where $0 \leq s \leq k{-}r{-}1$, we deduce that $\alpha = 0$.
\end{proof}

We next present a Khovanskii basis \cite{BPT} for the coordinate ring of the variety ${\rm SH}(k,n,r)$.
Khovanskii bases used to be called SAGBI bases in earlier works, and our
  arguments follows those given for Grassmannians in \cite[Section 3.1]{AIT} 
  and for flag varieties in \cite[Chapter 14]{CCA}.

We fix the reverse lexicographic term  order $>$ on the polynomial ring $S = \CC[{\bf x}]$,
where $x_{11} > x_{12} > \cdots > x_{n-k+r,n}$.  This is a {\em diagonal term order}, i.e.~for each minor of ${\bf x}$, the initial monomial is the product of the  entries on its diagonal.
Our coordinate ring is the image of the polynomial ring $R = \CC\bigl[ \langle I \rangle, [ J ]\, \bigr]$
under the ring homomorphism $\phi = \phi_{k,n,r}$ into $S$. 
  See Remark \ref{rmk:kernel}.
 For each of the $2 \binom{n}{k}$ generators of $R$,
we consider the initial monomial of its image in $S$. This gives a list of
$\binom{n}{k}$ monomials $ {\rm in}_> \phi( \langle I \rangle ) $ of degree $k$
and $\binom{n}{k}$ monomials  ${\rm in}_> \phi( [J ] ) $ of degree $n-k$.
These monomials lie in the initial algebra of our coordinate ring
\begin{equation}
\label{eq:initialalgebra} {\rm in}_> (\phi(S)) \,\, = \,\,
\CC \bigl[\, {\rm in}_>(f) \,: \, f \in \phi(S) \,\bigr]. 
\end{equation}

\begin{theorem} \label{thm:khovanskii}
The $2\binom{n}{k}$ minors $\phi(\langle I \rangle) $ and  $\phi([J])$ are a Khovanskii basis for the
coordinate ring $\phi(R)$ of the spinor-helicity variety ${\rm SH}(k,n,r)$, i.e.~their
initial monomials generate (\ref{eq:initialalgebra}).
\end{theorem}    

\begin{proof}    
Our argument mirrors that of \cite[Theorem 14.11]{CCA}.
We use the set-up in the proof of Theorem \ref{thm:main1}.
Monomials of bidegree $(d_1,d_2)$ in $S$ are represented by 
skew tableaux. These are formed by placing $d_1$ increasing rows
of length $n-k$, shifted by $r$ steps to the right,
above  $d_2$ increasing rows of length $k$.
This is seen in (\ref{eq:skewSymTableau}) for $d_1=d_2=1$.
Following \cite[Lemma 14.13]{CCA}, a monomial in $S$ is the
initial monomial of an element in $\phi(R)$
 if and only its representation as a skew tableau, as in \cite[eqn (14.4)]{CCA},
is a semi-standard skew tableau. Hence the
initial algebra (\ref{eq:initialalgebra}) is spanned as a
$\CC$-vector space by ${\bf x}$-monomials that correspond to
semi-standard skew tableaux. Every such monomial is
a product of diagonal monomials ${\rm in}_>\phi(\langle I \rangle) $  of degree $k$
and  diagonal  monomials ${\rm in}_>\phi([J])$ of degree $n-k$.
\end{proof}

We now illustrate Theorem \ref{thm:khovanskii} for the 
non-trivial instance shown in Figure~\ref{fig:P_(3,7,1)}.

\begin{example}[$k=3,n=7,r=1$]
The polynomial ring $S$ is generated by the entries~of 
\begin{equation}
\label{eq:35entries} 
{\bf x} \,\, = \,\,
\begin{pmatrix}
x_{11} & x_{12} &  x_{13} &  x_{14} &  x_{15} &  x_{16} &  x_{17} \\
x_{21} & x_{22} &  x_{23} &  x_{24} &  x_{25} &  x_{26} &  x_{27} \\
x_{31} & x_{32} &  x_{33} &  x_{34} &  x_{35} &  x_{36} &  x_{37} \\
x_{41} & x_{42} &  x_{43} &  x_{44} &  x_{45} &  x_{46} &  x_{47} \\
x_{51} & x_{52} &  x_{53} &  x_{54} &  x_{55} &  x_{56} &  x_{57} 
\end{pmatrix}.
\end{equation}
The polynomial ring $R$ is generated by the $70$ brackets in Figure~\ref{fig:P_(3,7,1)}.
The map $\phi: R \rightarrow S$~takes $\langle ij\ell \rangle$ to the $3 \times 3$-minor of ${\bf x}$ with row indices $\{1,2,3\}$ and column indices $\{i,j,\ell\}$. It takes $[ij\ell]$ to the signed $4 \times 4$-minor with row indices $\{2,3,4,5\}$
and column indices $[7] \backslash \{i,j,\ell\}$.
We consider the image (\ref{eq:initialalgebra}) of the map that takes each bracket to the diagonal initial monomial:
$$ \begin{matrix}
{\rm in}_>\phi \,:\, R \rightarrow S , & 
\!\!\!\! \langle 123 \rangle \,\mapsto \,x_{11} x_{22} x_{33}, &
\!\!\!\! \langle 124 \rangle \,\mapsto \,x_{11} x_{22} x_{34}, & \!\! \ldots\,, &
\!\!\!\! \langle 567 \rangle\, \mapsto \,x_{15} x_{26} x_{37}, \\ &
\, [ 123 ] \mapsto x_{24} x_{35} x_{46} x_{57} ,&
\, [ 124 ] \mapsto x_{23} x_{35} x_{46} x_{57} ,&  \!\! \ldots\,, &
\, [ 567 ] \mapsto x_{21} x_{32} x_{43} x_{54} .
\end{matrix}
$$
The kernel of the monomial map $\,{\rm in}_> \phi\,$ is a toric ideal in $R$.
This is minimally generated by
$329$ binomial quadrics. First, there are $140$ quadratic binomials from 
Young's poset $Y_{3,7}$:
\begin{equation}
\label{eq:140a}
\langle 125 \rangle \langle 134 \rangle - \langle 124 \rangle \langle 135 \rangle \, ,\,\,
 \langle 126 \rangle \langle 134 \rangle - \langle 124 \rangle \langle 136 \rangle \, ,\,\,
\ldots \, , \,
 \langle 367 \rangle \langle 457 \rangle - \langle 357 \rangle \langle 467 \rangle
\end{equation}
Likewise, the generators of
 the toric ideal  $\,{\rm kernel}\bigl({\rm in}_>(\phi)\bigr)\,$ include
  $140$  binomials  from~$\tilde Y_{3,7}$:
\begin{equation}
\label{eq:140b}
 [125] [134] - [124] [135]\,, \,\,
 [126] [134] - [124] [136]\,,\,  \ldots\, ,\,
 [367] [457] - [357] [467].
\end{equation}
Third, and most important, there are $49 $ mixed binomial quadrics in our toric ideal:
\begin{equation}
\label{eq:49mixed}
 \begin{matrix}
\langle 123 \rangle [123] - \langle 145 \rangle [145], & 
\langle 124 \rangle [123] + \langle 145 \rangle [135], & 
\langle 125 \rangle [123] - \langle 145 \rangle [134] , \\
\qquad \ldots \quad \ldots \quad \ldots \qquad,& 
\langle 123 \rangle [236] + \langle 124 \rangle [246], & 
\langle 123 \rangle [237] + \langle 124 \rangle [247] . \\
\end{matrix}
\end{equation}
The initial monomials in (\ref{eq:140a}), (\ref{eq:140b}) and (\ref{eq:49mixed})
are the incomparable pairs in the poset $\mathcal{P}_{3,7,1}$.
The toric variety defined by our binomials is a toric degeneration of
the spinor-helicity-variety ${\rm SH}(3,7,1)$. The
maximal chains of  $\mathcal{P}_{3,7,1}$ form a
triangulation of its {\em Newton-Okounkov body}.
\end{example}

We close this section with the remark that the poset $\Pcal_{k,n,r}$ is a distributive lattice, just like $Y_{k,n}$ and $\tilde Y_{k,n}$.
The join and meet operations $\wedge, \vee$ are defined as follows: if $\langle I \rangle [J]$ is an incomparable pair, then $\langle I  \vee J \rangle$ and $[I  \wedge J]$ are obtained by sorting the columns of the skew Young tableaux \eqref{eq:skewSymTableau}. 
With these lattice operations, the binomials (\ref{eq:49mixed}) 
can be written as
\[ 
    \langle I \rangle \cdot [ J ] \, - \, \langle I  \vee J \rangle \cdot [ I  \wedge J].
\]
In other words, the description in \cite[Theorem 14.16]{CCA} extends to the spinor-helicity varieties.

\section{Mandelstam Variety}
\label{sec4}

The componentwise multiplication of two vectors is known as the Hadamard product.
Applying this to pairs of
Pl\"ucker vectors gives rise to a rational map
\begin{equation}
\label{eq:hadamardmap}
 s \,:\,\PP^{\binom{n}{k}-1} \times \PP^{\binom{n}{k}-1} \,\dashrightarrow \,
\PP^{\binom{n}{k}-1}. 
\end{equation}
Generalizing the case $k=2$ in (\ref{eq:mandels}), the
coordinates of $s$ are called {\em Mandelstam invariants}:
\begin{equation}
\label{eq:mandelvar2}
 s_{i_1 i_2 \cdots i_k} \,\, = \,\, 
\langle   i_1 i_2 \cdots i_k \rangle  [ i_1 i_2 \cdots i_k ]. 
\end{equation}
We define the {\em Mandelstam variety} ${\rm M}(k,n,r)$  to be the closure of the image of the spinor-helicity variety 
${\rm SH}(k,n,r)$ under the Hadamard product map $s$.
Thus, ${\rm M}(k,n,r)$ is an irreducible variety in 
$\PP^{\binom{n}{k}-1}$.
We write  $\mathcal{I}({\rm M}(k,n,r))$ for the homogeneous prime ideal of
this variety. This
comprises all polynomial relations among the 
Mandelstam invariants $\, s_{i_1 i_2 \cdots i_k}$.

\begin{proposition} \label{prop:ambdim}
The linear span of the Mandelstam variety ${\rm M}(k,n,r)$ in $\PP^{\binom{n}{k}-1}$ 
is the subspace $\,\PP^N$  which is defined by the momentum conservation relations. Its dimension 
equals $$ N \,\,=\,\, \binom{n}{k}-1 - \binom{n}{k-r-1}. $$
\end{proposition}

This refers to the momentum conservation relations 
in the CEGM model \cite[eqn (5.6)]{CEGM}.

\begin{proof}
In our notation, the momentum conservation relations are written as follows:
\begin{equation}
\label{eq:momentumcons} \sum_{J \in \binom{[n]}{r+1}} \!\! s_{IJ} \,\, = \,\, 0
\qquad \hbox{for all} \,\,\, I \in \binom{[n]}{k-r-1}. 
\end{equation}
Note that, in the above sum, $s_{IJ} = 0$ whenever $I \cap J \neq \emptyset$. We claim that these linear forms  lie in
$\mathcal{I}({\rm M}(k,n,r))$ and that they are linearly independent.
To see this, recall the matrix $P Q^T$ from Theorem \ref{thm:PQ}.
The $\binom{n}{k-r-1}$ diagonal entries of~$PQ^T$ are 
$\, \sum_{J \in \binom{[n]}{r+1}} \langle IJ \rangle [IJ] $,
where the index $I$ runs over $ \binom{[n]}{k-r-1}$. 
This sum agrees with (\ref{eq:momentumcons}), which therefore lies
in $\mathcal{I}({\rm M}(k,n,r))$. The argument with the Johnson matrix in the
proof of Theorem  \ref{thm:PQ} shows that our $\binom{n}{k-r-1}$ linear forms
are linearly independent. The
dimension count in Corollary \ref{cor:arequadratic} implies that
they span the space of all linear forms in $\,\mathcal{I}({\rm M}(k,n,r))$.
\end{proof}

\begin{example}[$k=3,n=6,r=1$] \label{eq:M361a}
There are six momentum conservation relations:
\begin{equation}
\label{eq:MCR361}  \begin{matrix}
    s_{123}+s_{124}+s_{125}+s_{126}+s_{134}+s_{135}+s_{136}+s_{145}+s_{146}+s_{156}  &=& 0, \\
    s_{123}+s_{124}+s_{125}+s_{126}+s_{234}+s_{235}+s_{236}+s_{245}+s_{246}+s_{256} &=& 0, \\
    s_{123}+s_{134}+s_{135}+s_{136}+s_{234}+s_{235}+s_{236}+s_{345}+s_{346}+s_{356} &=& 0, \\
    s_{124}+s_{134}+s_{145}+s_{146}+s_{234}+s_{245}+s_{246}+s_{345}+s_{346}+s_{456} &=& 0, \\
    s_{125}+s_{135}+s_{145}+s_{156}+s_{235}+s_{245}+s_{256}+s_{345}+s_{356}+s_{456} &=& 0, \\
    s_{126}+s_{136}+s_{146}+s_{156}+s_{236}+s_{246}+s_{256}+s_{346}+s_{356}+s_{456} &=& 0.
  \end{matrix}
\end{equation}
  These define a subspace $\PP^{13}$ of $\PP^{19}$.
  The variety ${\rm M}(3,6,1)$ has codimension four in this $\PP^{13}$. 
A  general formula for the dimension of any Mandelstam variety is
given in the next result.
\end{example}

\begin{proposition} \label{prop:mandeldim}
The dimension of the Mandelstam variety equals
$$ {\rm dim}({\rm M}(k,n,r)) \,\, = \,\,
{\rm dim}({\rm SH}(k,n,r)) \, - n + 1 \,\, = \,\,
2k(n-k) - (k-r)^2 - n + 1 .$$
\end{proposition}

\begin{proof}
The $n$-dimensional torus $(\CC^*)^n$ acts on the spinor-helicity
variety as follows:
\begin{equation}
\label{eq:doubletorus}
 \begin{matrix} \langle i_1 i_2 \cdots i_k \rangle & \mapsto & 
t_{i_1} \, t_{i_2} \,\cdots \,t_{i_k} \,\langle i_1 i_2 \cdots i_k \rangle ,
\smallskip \\
[ i_1 i_2 \cdots i_k ] & \mapsto & 
t_{i_1}^{-1} t_{i_2}^{-1} \cdots t_{i_k}^{-1} [ i_1 i_2 \cdots i_k ].
\end{matrix} 
\end{equation}
The stabilizer is  a one-dimensional torus $\CC^*$.
The ring of polynomial invariants of the torus action is generated
by the Mandelstam invariants (\ref{eq:mandelvar2}).
Therefore, ${\rm M}(k,n,r)$ is the image of the quotient map, and its
dimension is $n-1$ less than the dimension of ${\rm SH}(k,n,r)$.
\end{proof}

Propositions  \ref{prop:ambdim} and \ref{prop:mandeldim}
are illustrated in Table \ref{tab:dim}. For the given values
of $k,n$ and $r$, we display the dimension
of ${\rm M}(k,n,r)$ and the dimension $N$
of its linear span. Note that $N$ can
be quite a bit smaller than the dimension
$\binom{n}{k}-1$ of the ambient Pl\"ucker space.
For example, ${\rm M}(3,8,0)$ has dimension $14$
inside a linear subspace $\PP^{27}$ of the Pl\"ucker space
$\PP^{55}$.

\begin{table}[h]
\setcounter{MaxMatrixCols}{20}
$$
\begin{matrix}
  & \!\! k,n=\!\!\!\! & 2, 4 &  2, 5 &  2, 6 &  2, 7 &  2, 8 &   3, 6 &  3, 7 &  3, 8 &  3, 9 & 4, 8 &  4, 9 
  \medskip \\
r=0 && 1, 1 &   4, 4  &   7, 8  &   10, 13  &   13, 19  &
   4, 4  &   9, 13  &   14, 27  &   19, 47  &   9, 13  &   16, 41  \\
r=1 &&   4, 4  &   7, 8  &   10, 13  &   13, 19  &   16, 26  &     9, 13  &   14, 27  &   19, 47  &   24, 74  &   16, 41  &   23, 89  \\
r=2 & &   5, 5  &   8, 9  &   11, 14  &   14, 20  &   17, 27  &    \! 12, 18  &   17, 33  &   22, 54  &   27, 82  &   21, 61  &   \!28, 116 
\end{matrix} \vspace{-0.3cm}
$$
\caption{\label{tab:dim}
The dimension of the Mandelstam variety ${\rm M}(k,n,r)$
and its ambient space $\PP^N$.}
\end{table}

We next discuss  the Mandelstam variety in the
case of primary interest in physics, namely $k=2$.
This lives in $\PP^{\binom{n}{2}-1}$.
 Here
$(s_{ij})$ is a symmetric $n \times n$ matrix with
 $s_{11} = \cdots = s_{nn} = 0$.

\begin{proposition}\label{prop:weak}
The $5 \times 5$ minors of $(s_{ij})$
vanish on the varieties ${\rm M}(2,n,r)$.
For $r=1$, the sum of all matrix entries is zero.
For $r=0$, each row and  each column in $(s_{ij})$ sums to zero, so we
only need the
$\frac{1}{2}\bigl(\binom{n-1}{5}^2+\binom{n-1}{5}\bigr)$ minors involving the last row or last column.
\end{proposition}

\begin{proof}
The sum constraints are the momentum conservations relations in  (\ref{eq:momentumcons}).
It suffices to prove the first sentence for $r=2$, when there are no such relations.
Note that ${\rm M}(2,n,2)$ is the Hadamard product \cite{BC} of ${\rm Gr}(2,n) $ with itself.
For each matrix $s $ in ${\rm M}(2,n,2)$, we have
\begin{equation}
\label{eq:rank4}
\begin{matrix} s_{ij} \,\,= \, \langle i\,j \rangle [ i \, j ] \,=\,
(\lambda_{1i} \lambda_{2j} -  \lambda_{2i} \lambda_{1j} )
(\tilde \lambda_{1i} \tilde \lambda_{2j} -  \tilde \lambda_{2i} \tilde \lambda_{1j} )
\qquad \qquad \qquad \qquad  \\ \qquad \qquad \quad
=\,  \lambda_{1i}  \tilde \lambda_{1i} \cdot \lambda_{2j} \tilde \lambda_{2j}
\,-\, \lambda_{1i} \tilde \lambda_{2i}  \cdot \lambda_{2j} \tilde \lambda_{1j}
\,-\,  \lambda_{2i}   \tilde \lambda_{1i} \cdot \lambda_{1j}  \tilde \lambda_{2j}
\,+\,  \lambda_{2i}  \tilde \lambda_{2i} \cdot \lambda_{1j} \tilde \lambda_{1j}.
\end{matrix}
\end{equation}
This shows that $s = (s_{ij})$ is a sum of four matrices of rank one, and hence
${\rm rank}(s) \leq 4$. 
\end{proof}

Our next theorem states that
the equations in Proposition \ref{prop:weak}
generate the prime ideals.

\begin{theorem} \label{thm:prime2n}
For $r=0,1,2$, the prime ideal of the Mandelstam variety ${\rm M}(2,n,r)$ is generated 
by the $5 \times 5$ minors of 
the matrix $(s_{ij})$ together with the 
respective linear
forms in~(\ref{eq:momentumcons}).
\end{theorem}

\begin{proof} The dimensions of our three Mandelstam varieties
from Proposition~\ref{prop:mandeldim} are
\begin{equation}
\label{eq:M2nrdim}
{\rm dim}({\rm M}(2,n,2)) = 3n-7 , \,\,
{\rm dim}({\rm M}(2,n,1)) = 3n-8 \quad {\rm and} \quad
{\rm dim}({\rm M}(2,n,0)) = 3n-11.
\end{equation}
Let $J$ denote the ideal generated  by the $5  \times 5$ minors.
This ideal is prime. This was shown for $4 \times 4$ minors in
\cite[Theorem 3.4]{DFRS}. The proof is the same for $5 \times 5$ minors.
A dimension count shows that $V(J)$ has dimension $3n-7$,
and from this we obtain the first assertion.

Every matrix in ${\rm M}(2,n,2)$ is a product $XX^T$
where  $X$ is an $n \times 4$ matrix whose rows   lie on the 
  Fermat quadric $\mathcal{V}(x_1^2 + x_2^2 + x_3^2 + x_4^2)$.
  Hence the coordinate ring of ${\rm M}(2,n,2)$ is the $n$-fold
  tensor product of the coordinate ring of the Fermat quadric.
  The latter ring is a normal domain and hence so is the former.
  All associated primes of the principal ideal
    generated by $\sum s_{ij}$ have height one in this domain.
  Using Lemma~\ref{lem:hope} below, we can now conclude that
this principal ideal is a prime ideal.
This argument shows that $J + \langle\, \sum s_{ij} \,\rangle$  
is a prime ideal in the polynomial ring $\CC[s]$.
The variety of this ideal has codimension $1$ in ${\rm M}(2,n,2)$.
Since this matches the dimension of ${\rm M}(2,n,1)$ 
in (\ref{eq:M2nrdim}), we conclude that $\mathcal{I}({\rm M}(2,n,1)) = 
J + \langle\, \sum s_{ij} \,\rangle$.

We now turn to $r=0$.
Let  $K$ be the ideal generated by $J$ and
$\sum_{j=1}^n s_{ij}$ for $i=1,2,\ldots,n$.
From Proposition \ref{prop:weak} we know that
${\rm M}(2,n,0) \subseteq \mathcal{V}(K)$.
We solve the $n$ linear equations for $s_{1n},s_{2n},\ldots,s_{n-1,n}$.
This leaves us with $\sum_{1 \leq i < j \leq n-1} s_{ij}=0$.
Moreover,
all $5 \times 5$ minors of the $n \times n$ matrix $s$ that involve the 
index $n$ are  sums of $5 \times 5$ minors that do not involve~$n$.
Using our result for $r=1$, we see that
$\CC[s]/K$ is isomorphic to the coordinate ring 
of ${\rm M}(2,n-1,1)$. This shows that $K$ is prime
and ${\rm dim}({\rm V}(K)) = 3(n-1)-8 = 3n-11$.
This matches the dimension for $r=0$ in
(\ref{eq:M2nrdim}). We thus conclude that $K$ is the prime ideal of ${\rm M}(2,n,0)$.
 \end{proof}

To complete the proof of Theorem \ref{thm:prime2n}, we still need to establish the following lemma.

\begin{lemma} \label{lem:hope}
The equation $\sum_{1 \leq i < j \leq n} s_{ij}  = 0$ defines a hypersurface that is reduced and irreducible in the variety of symmetric $n \times n$ matrices $(s_{ij})$ with zero diagonal and rank~$\leq 4$.
\end{lemma}

\begin{proof}
As in \cite[Theorem 3.4]{DFRS}, we work in the polynomial ring $\CC[X]$.
Our hypersurface is the variety cut out by the ideal $I$, which is  generated by
the quadric $\sum_{i,j=1}^n \sum_{k=1}^4 x_{ik} x_{jk}$ together with
the $n$ Fermat quadrics $\sum_{k=1}^4 x_{ik}^2$. 
These quadrics form a regular sequence in $\CC[X]$.
Hence $\CC[X]/I$ is a complete intersection ring.
By examining the Jacobian matrix of these $n+1$ quadrics, 
we can show that this variety is a complete intersection and that its
singular locus has codimension $\geq 2$. Serre's criterion 
for normality implies that the coordinate ring $\CC[X]/I$ is normal, so it
is a product of normal domains. This ring being graded, it has no non-trivial idempotents.
Hence it is a normal domain, so $I$ is a prime ideal. The argument just given is analogous to \cite[Theorem 2.1]{RSS}. We refer to that source for more details.
\end{proof}

We now turn to the general case $k  \geq 3$, with $r$ between $0$ and $k-2$.
The Mandelstam~invariants $s_{i_1 i_2 \cdots i_k}$ form a symmetric
$n \times  \cdots \times n$ tensor $s$, where an entry is zero whenever
$\# \{i_1,i_2 ,\ldots,i_k \} \leq k-1$. Its two-way marginal 
is the symmetric $n \times n$ matrix with entries 
\begin{equation}
\label{eq:matrixentry}  s_{i \,j \, +\cdots+} \,\,\,:= \,\,\,
\sum_{l_3=1}^n \sum_{l_4=1}^n \cdots \! \sum_{l_{k}=1}^n \! s_{i\,j\, l_3 l_4 \cdots l_k} 
\quad \hbox{for $\, 1 \leq i,j \leq n$.}
\end{equation}

\begin{proposition} \label{prop:marginal}
For every tensor $s$ in ${\rm M}(k,n,r)$, the
two-way marginal has rank~$\leq 4$.
\end{proposition}

\begin{example}[$k=3,n=6$] \label{eq:M361b}
The symmetric $6 \times 6 \times 6$ tensor $(s_{ijk}) $
has only $20$ distinct nonzero entries $s_{ijk}$ for $1\leq i < j< k \leq 6$. Its two-way marginal is the $6 \times 6$ matrix
\[ \!\!\!\!
\resizebox{\textwidth}{!}{$
\begin{bmatrix} 
0 \!\! & \!\! s_{123}{+}s_{124}{+}s_{125}{+}s_{126} \!\! & \!\! s_{123}{+}s_{134}{+}s_{135}{+}s_{136} \!\! & \!\! s_{124}{+}s_{134}{+}s_{145}{+}s_{146} \!\! & \!\! s_{125}{+}s_{135}{+}s_{145}{+}s_{156} \!\! & \!\! s_{126}{+}s_{136}{+}s_{146}{+}s_{156} \\
s_{123}{+}s_{124}{+}s_{125}{+}s_{126} \!\! & \!\! 0 \!\! & \!\! s_{123}{+}s_{234}{+}s_{235}{+}s_{236} \!\! & \!\! s_{124}{+}s_{234}{+}s_{245}{+}s_{246} \!\! & \!\! s_{125}{+}s_{235}{+}s_{245}{+}s_{256} \!\! & \!\! s_{126}{+}s_{236}{+}s_{246}{+}s_{256} \\
s_{123}{+}s_{134}{+}s_{135}{+}s_{136} \!\! & \!\! s_{123}{+}s_{234}{+}s_{235}{+}s_{236} \!\! & \!\! 0 \!\! & \!\! s_{134}{+}s_{234}{+}s_{345}{+}s_{346} \!\! & \!\! s_{135}{+}s_{235}{+}s_{345}{+}s_{356} \!\! & \!\! s_{136}{+}s_{236}{+}s_{346}{+}s_{356} \\
s_{124}{+}s_{134}{+}s_{145}{+}s_{146} \!\! & \!\! s_{124}{+}s_{234}{+}s_{245}{+}s_{246} \!\! & \!\! s_{134}{+}s_{234}{+}s_{345}{+}s_{346} \!\! & \!\! 0 \!\! & \!\! s_{145}{+}s_{245}{+}s_{345}{+}s_{456} \!\! & \!\! s_{146}{+}s_{246}{+}s_{346}{+}s_{456} \\
s_{125}{+}s_{135}{+}s_{145}{+}s_{156} \!\! & \!\! s_{125}{+}s_{235}{+}s_{245}{+}s_{256} \!\! & \!\! s_{135}{+}s_{235}{+}s_{345}{+}s_{356} \!\! & \!\! s_{145}{+}s_{245}{+}s_{345}{+}s_{456} \!\! & \!\! 0 \!\! & \!\! s_{156}{+}s_{256}{+}s_{356}{+}s_{456} \\
s_{126}{+}s_{136}{+}s_{146}{+}s_{156} \!\! & \!\! s_{126}{+}s_{236}{+}s_{246}{+}s_{256} \!\! & \!\! s_{136}{+}s_{236}{+}s_{346}{+}s_{356} \!\! & \!\! s_{146}{+}s_{246}{+}s_{346}{+}s_{456} \!\! & \!\! s_{156}{+}s_{256}{+}s_{356}{+}s_{456} \!\! & \!\! 0
\end{bmatrix}\!\!.$}
\]
This matrix has rank four on ${\rm M}(3,6,1)$ and hence also on ${\rm M}(3,6,0)$,
but not on ${\rm M}(3,6,2)$.
Geometrically, this matrix encodes the map
${\rm M}(3,6,1) \dashrightarrow {\rm M}(2,6,0)$ of Mandelstam varieties.
\end{example}

The Mandelstam variety ${\rm M}(k,n,r)$ has a natural parametrization,
namely the composition of the Hadamard map 
with the parametrization given by the  map $\phi_{k,n,r}$
in Remark~\ref{rmk:kernel}.
The parameters are the entries
in the $(n-k+r)  \times n$ matrix ${\bf x} = (x_{ij})$.
The Mandelstam invariant $s_I$ is a polynomial in the entries of ${\bf x}$.
Namely, $s_I$ is the product of
the  $k \times k$ minor indexed by $I$ in the first $k$ rows of ${\bf x}$
with the signed $(n - k) \times (n-k)$ minor indexed by $[n] \backslash I$
in the last $n-k$ rows of ${\bf x}$. The corresponding ring map $\phi_{k,n,r} \circ s^* $
has kernel $\mathcal{I} ({\rm M}(k,n,r))$.

\begin{proof}[Proof of Proposition \ref{prop:marginal}]
Since ${\rm M}(k,n,r) \subset {\rm M}(k,n,k-2)$ for all $r < k-2$, it suffices to prove the statement 
in the case $r=k-2$. For  this case, we consider the rational map
\begin{equation}
\label{eq:k-2map}
\SH(k,n,k-2) \,\dashrightarrow \,\SH(2,n,0) \,, \quad
(V,W) \,\mapsto \,(V \cap W^\perp, V^\perp \cap W).
\end{equation}
Let $\textbf{x}$ be the $(n-2) \times n$ matrix in Remark \ref{rmk:kernel} used to parametrize $\SH(k,n,k-2)$. In terms of the parametrization in Remark \ref{rmk:kernel}, the space
$V^\perp \cap W = (V + W^\perp)^\perp$ is the kernel of the matrix $X$, while
$V\cap W^\perp$ is the span of the $2$ middle rows of $\textbf{x}$ indexed by $\{k-1,k\}$.
So the Pl\"ucker coordinates of $V \cap W^\perp$ are exactly the $2 \times 2$ minors of the two middle rows, while the coordinates of $V^\perp \cap W$ can be obtained as signed maximal minors of the matrix $\textbf{x}$. Explicitly, the $ij$-th Pl\"ucker coordinate of $V^\perp \cap W$ is the signed minor of $\textbf{x}$ indexed by $[n]\setminus ij$. The map (\ref{eq:k-2map}) is equivariant with respect to the involution in Remark \ref{rmk:involution}, and it will induce a rational map of Mandelstam~varieties
\begin{equation}
\label{eq:k-2mapb}
{\rm M}(k,n,k-2)\, \dashrightarrow \,{\rm M}(2,n,0), \quad (s_{l_1 \cdots l_k}) \mapsto \left(\sum_{1 \leq l_3, \dots, l_k \leq n} s_{ijl_3l_4 \dots l_k} \right)_{1 \leq i,j \leq n}.
\end{equation}

The fact that the map (\ref{eq:k-2mapb}) is well defined and compatible with (\ref{eq:k-2map}) is non-trivial. To show this, we first consider the $(n-2) \times (2n-4)$ matrix $\tilde{\bf x}$ obtained by horizontally concatenating the two matrices ${\bf x}$ and
\[
\begin{bmatrix}
{\rm id}_{k-2}                      &   0_{(k-2) \times (n-k-2)}   \\    
   0_{2 \times (k-2)}              &    0_{2 \times (n-k-2)} \\    
   0_{(n-k-2) \times (k-2)}   &   {\rm id}_{n-k-2}  \\        
\end{bmatrix}.
\]
We then note that the Pl\"ucker coordinates of the spaces $V, W, V \cap W^\perp$ and $V^\perp \cap W$ can be written as maximal minors of $\tilde{\bf x}$, i.e. as coordinates on $\Gr(n-2, 2n-4)$. Writing everything in terms of these minors, and using the straightening laws on the Grassmannian $\Gr(n-2,2n-4)$, we can show that, when $(V,W) \in \SH(k,n,k-2)$, the right hand side of (\ref{eq:matrixentry}) factors and is equal to $s_{ij} = \langle ij \rangle [ij]$ for the pair $(V\cap W^\perp, V^\perp \cap W) \in \SH(2,n,0)$.
For example, when $n=6,k=3,i=1,j=2$, the desired identity $\, s_{12} \,= \, s_{123}  + s_{124} + s_{125} + s_{126} \,$ is precisely the five-term relation that holds on $\Gr(4,8)$ in Example \ref{ex:foureight}. This proves that the map \eqref{eq:k-2mapb} is well defined and compatible with the map \eqref{eq:k-2map}. We then conclude the proof by noting that the matrix $(s_{ij})$ has rank $\leq 4$, by virtue of Proposition~\ref{prop:weak}.
\end{proof}

\begin{remark}
    The rational map in \eqref{eq:k-2map} is a well-defined morphism on the 
    open subset $\SH(k,n,k-2) \backslash \SH(k,n,k-3)$. This open subset is 
     the smooth locus of $\SH(k,n,k-2)$. In general, the singular locus of
     the spinor-helicity variety $\SH(k,n,r)$ is precisely $\SH(k,n,r-1)$.
\end{remark}

 We now replace $\phi_{k,n,r}$ with a birational parametrization
 for (\ref{eq:rephrasing}), namely
   by specializing~${\bf x}$ as follows.
  Start  rows $r+1,\ldots,k$ with a unit matrix.
This leaves $(k-r)(n-k+r)$ parameters for $V \cap W^\perp$.
Start rows $1,\ldots,r$ with $k-r$ zero columns,
followed by a unit matrix. This leaves $ r (n-k)$
parameters for $V/(V \cap W^\perp)$.
Start rows $k+1,\ldots,n-k+r$
with $k-r$ zero columns, followed by a unit matrix.
This leaves $(n-2k+r) k $
parameters for $W^\perp/(V \cap W^\perp)$.

\begin{corollary} \label{cor:rulesabove}
The rules above give a birational map 
$\psi_{k,n,r} : \CC^{2k(n-k) -(k-r)^2} \! \to {\rm SH}(k,n,r).$
\end{corollary}

\begin{proof} 
First note that the parameter count above matches the dimension formula in 
(\ref{eq:dimformula}):
$$
(k-r)(n-k+r) \, + \,
r(n-k) \, + \,
(n-2k+r) k  \,\,=\,\,  2k(n-k) - (k-r)^2.$$
To parametrize (\ref{eq:rephrasing}), one first
chooses $V \cap W^\perp$, and thereafter
$V/(V \cap W^\perp)$ and $W^\perp/(V \cap W^\perp)$.
Each block of rows gives a
birational parametrization of the respective Grassmannian.
\end{proof}

By composing $\psi_{k,n,r}$ with the Hadamard map $s$, we
obtain a parametrization of the Mandelstam variety ${\rm M}(k,n,r)$.
Each fiber has dimension $n-1$, reflecting the torus action
in (\ref{eq:doubletorus}).
To obtain a finite-to-one parametrization of
${\rm M}(k,n,r)$,
we can now replace $n-1$ of the matrix
entries $x_{ij}$ by $1$.
The resulting parametric representation of  ${\rm M}(k,n,r)$
can then be used in
numerical algebraic geometry. The following proposition serves
as an illustration.

\begin{proposition}
\label{prop:M361}
The Mandelstam variety ${\rm M}(3,6,1)$
has dimension $9$ and degree $56$ in~$\PP^{19}$.
Its prime ideal is minimally generated by $14$ quartics,
plus the six linear forms in (\ref{eq:MCR361}).
Inside their subspace $\PP^{13}$,
the variety is arithmetically Cohen-Macaulay, and its Betti diagram equals
 \begin{equation}
 \label{eq:bettidiagram} \, 
 \begin{bmatrix} \,14 & \cdot & \cdot & \cdot \,\\
\, \cdot & 56 & 64 & 21\, 
\end{bmatrix}. 
\end{equation}
 \end{proposition}
 
 \begin{proof}[Computational proof]
 The map $\psi_{3,6,1}$ is given by
     the following specialization of our matrix:
 $$ {\bf x} \,\, = \,\, \begin{bmatrix}
\,0 & 0 & 1 & b_1 & b_2 & b_3 \\
\,1 & 0 & a_1 & a_2 & a_3 & a_4 \\
\,0 & 1 & a_5 & a_6 & a_7 & a_8 \\
\,0 & 0 & 1 & c_1 & c_2 & c_3 \\
\end{bmatrix}.
$$
We now remove five parameters by setting $a_1=a_5 = c_1 = c_2 = c_3 = 1$.
The Mandelstam invariants are polynomials in the nine remaining unknowns 
$a_i$ and $b_j$. This specifies a
two-to-one map
$\,\CC^9 \rightarrow \PP^{19}$.
A computation checks that the Jacobian of the map has full rank.
Hence the closure of its image is the $9$-dimensional Mandelstam variety ${\rm M}(3,6,1)$.
This lies  in the $\PP^{13}$ defined by (\ref{eq:MCR361}). 
We use these relations to eliminate the six variables
$\,s_{123},s_{124}, s_{134},s_{234},s_{235},s_{236}$.
Thereafter, we can
 view $ {\mathcal I}({\rm M}(3,6,1))$ as an ideal
in the remaining $14$ variables. 
Our computations with this ideal were carried out in
 {\tt Macaulay2}~\cite{M2}.

The ideal contains no quadrics or cubics, but we find $14$
 linearly independent quartics. Let $I$ be the subideal
 generated by the $14$ quartics. This  has codimension 
 $4$ and degree $56$, and it
contains the $21$ quintics given by the $5 \times 5$ minors of the
$6 \times 6$ matrix in Example~\ref{eq:M361b}.
We compute the minimal free resolution of $I$, and find that
its Betti diagram equals (\ref{eq:bettidiagram}).
Thus, $I$ is Cohen-Macaulay, and so it
is an intersection of primary ideals of codimension $4$.

We now apply {\tt HomotopyContinuation.jl} \cite{BT} to the map
$\,\CC^9 \rightarrow \PP^{19}$,
and we compute the degree of its image.
This yields an independent proof that   ${\rm M}(3,6,1) $ has degree $56$.
Since $I$ has degree $56$,  and since the degree is additive over
the primary components,
this shows that the ideal $I$ is prime.
We conclude that $I =  {\mathcal I}({\rm M}(3,6,1))$, and the proof is complete.
\end{proof}

\begin{remark}
At present, we do not know the meaning of our $14$ quartic generators.
 The shortest quartic we found has $140$ monomials.
In reverse lexicographic order, it looks like
\[
\resizebox{\textwidth}{!}{$
\begin{matrix}
s_{136} s_{156} s_{235} s_{345}-s_{135} s_{156} s_{236} s_{345}-s_{145} s_{156} s_{236} s_{345}-s_{156}^2 s_{236} s_{345}-s_{146} s_{156} s_{245} s_{345} \phantom{moo}
 \\ +\,s_{135} s_{156} s_{246} s_{345} 
+s_{145} s_{156} s_{246} s_{345}+s_{156}^2 s_{246} s_{345}+s_{136} s_{145} s_{256} s_{345}-s_{135} s_{146} s_{256} s_{345}+ \,\cdots 
\\
\cdots \,\, \cdots \,\,+\,s_{136} s_{345} s_{456}^2+s_{156} s_{345} s_{456}^2-s_{135} s_{346} s_{456}^2-2 s_{135} s_{356} s_{456}^2+s_{145} s_{356} s_{456}^2-s_{135} s_{456}^3.
\end{matrix}$}
\]
\end{remark}

We close this section by recording a few more general facts about Mandelstam varieties.

\begin{proposition} \label{prop:deletion}
The Mandelstam variety ${\rm M}(k,n,k)$ is the Hadamard square
of the Grassmannian ${\rm Gr}(k,n)$. It contains all other Mandelstam 
varieties by the chain of inclusions
\begin{equation}
\label{eq:endsec4a}
{\rm M}(k,n,0) \,\subset\, {\rm M}(k,n,1) \,\subset \,
{\rm M}(k,n,2) \,\subset\, \cdots \,\subset\, {\rm M}(k,n,k).
\end{equation}
There is natural chain of 
dominant deletion maps, induced by removing columns in $\lambda$ and $\tilde \lambda$:
\begin{equation}
\label{eq:endsec4b}
{\rm M}(k,n,0) \,\dashrightarrow \, {\rm M}(k,n-1,1) \,\dashrightarrow \,
{\rm M}(k,n-2,2) \,\dashrightarrow \,\, \cdots \,\,\dashrightarrow \, {\rm M}(k,n-k,k).
\end{equation}
\end{proposition}

\begin{proof}[Proof and discussion]
The term {\em Hadamard square} refers to the Hadamard product of a variety with itself.
For an introduction to Hadamard products of varieties see the book \cite{BC}.
We obtain inclusions ${\rm SH}(k,n,r) \subset {\rm SH}(k,n,r+1)$
by relaxing the rank constraints in (\ref{eq:rephrasing}),
and we obtain surjections 
${\rm SH}(k,n,r) \rightarrow {\rm SH}(k,n-1,r+1)$
by deleting the last columns in $\lambda$ and $\tilde \lambda$ respectively.
These maps are compatible with Hadamard products, so they descend
to inclusions 
${\rm M}(k,n,r) \subset {\rm M}(k,n,r+1)$
and surjections
${\rm M}(k,n,r) \rightarrow {\rm M}(k,n-1,r+1)$.
It would be interesting to study the fibers of the maps in (\ref{eq:endsec4b}).
Their dimensions are $0,2,4,\ldots,2k-2$.
\end{proof}

\section{Positivity and Tropicalization}
\label{sec5}

In our last two sections, we set the stage for future research
on spinor-helicity varieties.
Our view is now aimed towards
tropical geometry, positive geometry, and applications to scattering~amplitudes.

Bossinger, Drummond and Glew \cite{BDG} studied
the Gr\"obner fan and positive geometry of the variety ${\rm SH}(2,5,0)$
in Example \ref{ex:twofivezero}
which they identified with the
Grassmannian ${\rm Gr}(3,6)$.
We shall examine this in a broader context.
The following result explains their identification.

\begin{proposition} \label{prop:larageneral}
    For any $k \geq 1$, 
    the varieties $\SH(k,2k{+}1,0)$ and $\SH(k{+}1,2k{+}1,1)$ are isomorphic and their coordinate ring is isomorphic to that of the Grassmannian $\Gr(k{+}1, 2k{+}2)$.
    \end{proposition}
    
\begin{proof}
The isomorphism between $\SH(k, 2k\!+\!1,0)$ and $\SH(k\!+\!1, 2k\!+\!1,1)$ 
arises because $(V,W) $ is in $\SH(k\!+\!1, 2k\!+\!1,1)$ if and only if $(V^\perp, W^\perp) 
$ is in $ \SH(k,2k\!+\!1,0)$.
The identification with the Grassmannian
${\rm Gr}(k\!+\!1,2k\!+\!2)$ uses the specialized parametrization $\psi_{k,2k+1,0}$
in Corollary \ref{cor:rulesabove}. We introduce a new parameter $z$,
to account for the fact that 
the dimension of ${\rm Gr}(k\!+\!1, 2k\!+\! 2)$ exceeds the dimension of ${\rm SH}(k,2k\!+\!1) $ by one. We also
augment ${\bf x}$  with one extra column $(0,0,\ldots,0,z)^T$.
The new matrix ${\bf x}$ has $k+1$ rows and $2k+2$ columns,
and it contains $(k+1)^2$ parameters. 
Its maximal minors give a
 birational parametrization of ${\rm Gr}(k\!+\!1, 2k\!+\! 2)$
and also of  $\SH(k, 2k\!+\!1,0)$.
The minors involving the extra column are the
Pl\"ucker coordinates for $V$. The others
are  Pl\"ucker coordinates for $W^\perp$.
\end{proof}

The {\em positive Grassmannian} ${\rm Gr}_+(k,n)$
is defined by requiring all Pl\"ucker coordinates
$\langle i_1 i_2 \cdots i_k \rangle$ of the subspace $V$ to be real and positive.
We define the {\em dually positive Grassmannian} ${\rm Gr}^+(k,n)$
to be  ${\rm Gr}_+(n-k,n)$ under the identification between $W$ and $W^\perp$.
Thus ${\rm Gr}^+(k,n)$ is an open semialgebraic set 
isomorphic to ${\rm Gr}_+(n-k,n)$. It is defined by
\begin{equation}
\label{eq:Grplus}
{\rm sign}\bigl( [j_1 j_2 \cdots j_k] \bigr) \,= \,(-1)^{j_1+j_2 + \cdots + j_k}
\quad \hbox{for} \quad
1 \leq j_1 < j_2 < \cdots < j_k \leq n.
\end{equation}
The {\em positive spinor-helicity variety} consists of all positive points in 
our variety:
\begin{equation}
\label{eq:SHplus} {\rm SH}_+(k,n,r) \,\,\, := \,\,\,
{\rm SH}(k,n,r) \,\, \cap \,\, \bigl(\,{\rm Gr}_+(k,n) \,\times\, {\rm Gr}^+(k,n) \bigr)
\quad \subset \,\,\,
 \RR \PP^{\binom{n}{k}-1} \times 
 \RR \PP^{\binom{n}{k}-1} .
\end{equation}
We finally define the {\em positive Mandelstam variety}
$\,{\rm M}_+(k,n,r)$  by the inequalities in (\ref{eq:Grplus}):
\begin{equation}
\label{eq:Mplus}
{\rm sign}\bigl( s_{j_1 j_2 \cdots j_k} \bigr) \,= \,(-1)^{j_1+j_2 + \cdots + j_k}
\quad \hbox{for} \quad
1 \leq j_1 < j_2 < \cdots < j_k \leq n .
\end{equation}
Thus ${\rm M}_+(k,n,r)$ is a semialgebraic subset of $\PP^{\binom{n}{k}-1}$.
It contains the image of ${\rm SH}_+(k,n,r)$ under the 
Hadamard product map $s$ in (\ref{eq:hadamardmap}). In general, this inclusion is strict. For example,
${\rm M}_+(2,4.2)$ strictly contains the Hadamard product of $\Gr_+(2,4)$ and $\Gr^+(2,4)$. To see this, we note that the following expression is positive on the latter set but not on the former set:
\begin{equation}
\label{eq:isstrict}
  s_{13} s_{24} + s_{14} s_{23} - s_{12} s_{34}  
  \,\, = \,\,    \langle13\rangle \langle24\rangle [14] [23] \,+\, \langle14\rangle \langle23\rangle [13] [24].
\end{equation}
  
We now recycle Example~\ref{ex:twofivezero} and Proposition~\ref{prop:weak}
for our running example in this section.

\begin{example}[$k\!=\!2,n\!=\!5$] \label{ex:running5a}
For points in ${\rm SH}_+(2,5,r)$, the
rank $2$ matrices $P$ and $Q$  satisfy
\[
{\rm sign}(P) \,\, =\,\,
\begin{footnotesize}
\begin{bmatrix} 
0 & + & + & + & + \\
 - & 0 &  + & + & + \\
 - & - & 0 & + & + \\
 - & - & - & 0 & + \\
 - & - & - & - & 0
\end{bmatrix}
\end{footnotesize}
\qquad {\rm and} \qquad
{\rm sign}(Q) \,\, =\,\,
\begin{footnotesize}
 \begin{bmatrix} 
0 & - & + & - & + \\
+ & 0 &  - & + & - \\
 - & + & 0 & - & + \\
 + & - & + & 0 & - \\
 - & + & - & + & 0
\end{bmatrix}.
\end{footnotesize}
\]
The positive spinor-helicity variety ${\rm SH}_+(2,5,0)$ is the subset defined by the
 equations in (\ref{eq:momentum}).

 The positive Mandelstam variety
${\rm M}_+(2,5,2) $ is a $9$-dimensional simplex $\RR \PP_+^9$.
It consists of 
symmetric $5 \times 5$ matrices $s = (s_{ij})$ with alternating sign pattern,
i.e.  ${\rm sign}(s) = {\rm sign}(Q)$. 
Its subset ${\rm M}_+(2,5,0)$ is defined in this simplex 
by the equations $\sum_{j=1}^5  s_{ij}  = 0$ for $i=1,2,3,4,5$.
This is a cyclic $4$-polytope with $6$ vertices. It has the
 f-vector $(6, 15, 18, 9)$. The vertices~are:
 
 \[
\resizebox{\textwidth}{!}{$        
\begin{pmatrix} 
        0 & - & + & 0 & 0 \\
        - & 0 &  0 & + & 0 \\
        + & 0 & 0 & - & 0 \\
        0 & + & - & 0 & 0 \\
        0 & 0 & 0 & 0 & 0
        \end{pmatrix},
        \begin{pmatrix} 
        0 & 0 & + & - & 0 \\
        0 & 0 &  - & + & 0 \\
        + & - & 0 & 0 & 0 \\
        - & + & - & 0 & 0 \\
        0 & 0 & 0 & 0 & 0
        \end{pmatrix},
        \begin{pmatrix} 
        0 & - & 0 & 0 & + \\
        - & 0 & 0 & + & 0 \\
        0 & 0 & 0 & 0 & 0 \\
        0 & + & 0 & 0 & - \\
        + & 0 & 0 & - & 0
        \end{pmatrix},
        \begin{pmatrix} 
        0 & 0 & 0 & - & + \\
        0 & 0 & 0 & + & - \\
        0 & 0 & 0 & 0 & 0 \\
        - & + & 0 & 0 & 0 \\
        + & - & 0 & 0 & 0
        \end{pmatrix},
         \begin{pmatrix} 
        0 & 0 & 0 & 0 & 0 \\
        0 & 0 & - & + & 0 \\
        0 & - & 0 & 0 & + \\
        0 & + & 0 & 0 & - \\
        0 & 0 & + & - & 0
        \end{pmatrix},
        \begin{pmatrix} 
        0 & - & + & 0 & 0 \\
        - & 0 & 0 & + & 0 \\
        + & 0 & 0 & - & 0 \\
        0 & + & - & 0 & 0 \\
        0 & 0 & 0 & 0 & 0
        \end{pmatrix}.$}
\]
Nine of the inequalities ${\rm sign}(s_{ij}) = (-1)^{i+j}$ define facets.
Only $s_{24} \geq 0$ is not facet-defining.
\end{example}

The second thread in this section is {\em tropical geometry}.
Using notation from the textbook \cite{MStrop},
the tropicalizations of our varieties ${\rm SH}(k,n,r)$ 
and ${\rm M}(k,n,r)$ are the tropical varieties
\begin{equation}
\label{eq:tropvar}
{\rm trop}({\rm SH}(k,n,r)) \,\,\subset\,\,
\RR^{\binom{n}{k}}\! / \RR{\bf 1} \,\times \,\RR^{\binom{n}{k}} \! / \RR {\bf 1}
\,\,\quad {\rm and} \quad\,\,
{\rm trop}({\rm M}(k,n,r)) \,\,\subset \,\,
 \RR^{\binom{n}{k}}\! / \RR {\bf 1}.
\end{equation}
 These are balanced polyhedral fans  whose dimensions
 are given by  Propositions  \ref{prop:geoflag} and~\ref{prop:mandeldim}.
 Each such fan is a finite intersection of the tropical hypersurfaces given by
a tropical basis.
We illustrate these concepts for an example where the
underlying variety is a linear space.

\begin{example}[$k=2,n=5$] \label{ex:running5b} 
The tropical linear space ${\rm trop}({\rm M}(2,5,0))$ is a 
pointed fan of dimension $4$ in $\RR^{10}/\RR {\bf 1}$.
A minimal tropical basis consists of $15$ linear forms with four terms:
\[
\resizebox{\textwidth}{!}{$
\begin{matrix}
s_{12}{+}s_{13}{+}s_{14}{+}s_{15},
 s_{12}{+}s_{23}{+}s_{24}{+}s_{25},
 s_{13}{+}s_{23}{+}s_{34}{+}s_{35},
 s_{14}{+}s_{24}{+}s_{34}{+}s_{45},
 s_{15}{+}s_{25}{+}s_{35}{+}s_{45}, \\
s_{12}{+}s_{13}{+}s_{23}{-}s_{45},
s_{12}{+}s_{14}{+}s_{24}{-}s_{35},
s_{12}{+}s_{15}{+}s_{25}{-}s_{34},
s_{34}{+}s_{35}{+}s_{45}{-}s_{12},
s_{13}{+}s_{14}{+}s_{34}{-}s_{25}, \\
s_{13}{+}s_{15}{+}s_{35}{-}s_{24},
s_{24}{+}s_{25}{+}s_{45}{-}s_{13},
s_{14}{+}s_{15}{+}s_{45}{-}s_{23},
s_{23}{+}s_{25}{+}s_{35}{-}s_{14},
s_{23}{+}s_{24}{+}s_{34}{-}s_{15}.
\end{matrix}$}
\]
The underlying rank $5$ matroid on $10$ elements is  the
exceptional unimodular matroid~$R_{10}$.
This matroid appears in \cite[Section 3.3]{YY}.
Hence ${\rm trop}({\rm M}(2,5,0))$ is the cone
over the Bergman complex of $R_{10}$. This complex consists of $315$ tetrahedra and $45$ bipyramids and its f-vector is $(40,240,510,360)$. 
The $40$ vertices are the $10$ coordinate points $e_{ij}$ and the $30$ circuits of~$R_{10}$.
\end{example}

Combining our two threads, and using the notion of
positivity defined above leads us to
\begin{equation}
\label{eq:tropvarplus}
{\rm trop}_+({\rm SH}(k,n,r)) \,\,\subset\,\,
\RR^{\binom{n}{k}}\! / \RR {\bf 1} \,\times \,\RR^{\binom{n}{k}} \! / \RR {\bf 1}
\,\,\quad {\rm and} \quad\,\,
{\rm trop}_+({\rm M}(k,n,r)) \,\,\subset \,\,  \RR^{\binom{n}{k}}\! / \RR {\bf 1}.
 \end{equation}
 These positive tropical varieties are subfans of the respective tropical varieties.
 They are defined by requiring sign compatibility in the tropical equations.
 For details see \cite{AKW, BLS, SW}.
 
\begin{example}[$k=2,n=5$]  \label{ex:running5c}
Following Examples \ref{ex:running5a} and \ref{ex:running5b},
the positive tropical Mandelstam variety ${\rm trop}_+({\rm M}(2,5,0))$ is 
defined by the following system of $15$ tropical equations:
\[
\begin{matrix}
s_{13} \oplus s_{15} = s_{12} \oplus s_{14}\,,\,\,
s_{24} = s_{12} \oplus s_{23} \oplus s_{25}\,,\,\,
s_{13} \oplus s_{35} = s_{23} \oplus s_{34}\,,\,\,
s_{24} = s_{14} \oplus s_{34} \oplus s_{45}, \\
s_{15} \oplus s_{35} = s_{25} \oplus s_{45}\,, \,\,
s_{13} \oplus s_{45} = s_{12} \oplus s_{23}\,,\,\,
s_{24} = s_{12} \oplus s_{14} \oplus s_{35}\,,\,\,
s_{15} \oplus s_{34} = s_{12} \oplus s_{25}, \\
s_{12} \oplus s_{35} = s_{34} \oplus s_{45}\,,\,\,
s_{13} \oplus s_{25} = s_{14} \oplus s_{34}\,, \,\,
s_{13} \oplus s_{15} \oplus s_{35} = s_{24}\,,\,\,
s_{24} = s_{13} \oplus s_{25} \oplus s_{45}, \\
s_{15} \oplus s_{23}\, =\, s_{14} \oplus s_{45}\,,\,\,\,
s_{14} \oplus s_{35} \,=\, s_{23} \oplus s_{25}\,, \,\,\,
s_{24} \,=\, s_{15} \oplus s_{23} \oplus s_{34},
\end{matrix}
\]
Here $\,x \oplus y := {\rm min}(x,y)$.
Note the special role of the non-facet variable $s_{24}$.
These $15$ equations define the positive Bergman complex,
in the sense  of Ardila, Klivans and Williams~\cite{AKW}.

We find that $\,{\rm trop}_+({\rm M}(2,5,0))$ is the cone over
 a $3$-sphere. That sphere is glued from
$48$ tetrahedra and $18$ bipyramids, and its f-vector is $(24,108,150, 66)$.
Geometrically, this is a subdivision of  the boundary of
the $4$-polytope $\Delta_2 \times \Delta_2$, which is the product of
two triangles. This is dual to the cyclic polytope 
in Example \ref{ex:running5a}, and its
 f-vector is $(9,18,15,6)$.
The ``fine subdivision'' discussed in \cite[Corollary 3.5]{AKW}
is the barycentric subdivision of  ${\rm Bdr}(\Delta_2 \times \Delta_2)$.
\end{example}
 
\begin{remark}
The positive and tropical geometry of Grassmannians and
flag varieties has been studied intensely in recent years.
See \cite{Bor, BEZ, Ola, SW} for some references, and
\cite{CEGM, DFGK, Lam2} for physics perspectives.
We know from Proposition \ref{prop:geoflag} that
${\rm SH}(k,n,0) = {\rm Fl}(k,n-k;\CC^n)$ and
${\rm SH}(k,n,k) = {\rm Gr}(k,n) \times {\rm Gr}(k,n)$.
The positive geometry structure on partial flag varieties is well established in the literature; see, for example, \cite[Section 3.4]{Lam1}.
However, the notion of positivity for which the boundary structure and canonical form of partial flag varieties has been studied is the notion of \emph{total positivity} due to Lusztig.
On the other hand, the positive region $\SH_+(k,n,r)$ is the region where all the Pl\"ucker coordinates are nonnegative.
By virtue of \cite[Theorem 1.2]{Bloch_Karp}, these two notions of positivity on ${\rm Fl}(k,n-k; \CC^n)$ agree if and only if the dimensions $k$ and $n-k$ are equal or consecutive, i.e., $n = 2k$ or $n = 2k + 1$. Hence, at present, the positive region $\SH_+(k,n,r)$ is well understood only in the cases where $r = k$, or $r = 0$ with $n = 2k$ or $n = 2k + 1$.
Determining the boundary structure and canonical form of $\SH_+(k,n,r)$ for general $k,n,r$ remains an open problem, as does the study of ${\rm trop}_+({\rm SH}(k,n,r))$. Finally, we note that a detailed study of ${\rm trop}_+ ({\rm SH}(1,n,0))$ was carried out by Olarte in \cite{Ola}.
It would be desirable to extend this to $k \geq 2$ using the techniques introduced in \cite{Bor, BEZ}.
\end{remark}


\begin{corollary}  \label{cor:C2}
Modulo a scaling action by $\RR^+$, the positive Grassmannian $\Gr(k{+}1, 2k{+}2)$ coincides with the positive spinor-helicity varieties $\SH_+(k,2k{+}1,0)$ and $\SH_+(k{+}1,2k{+}1,1)$, for all $k \geq 1$. The analogous statement holds for their (positive) tropical varieties.
\end{corollary}
 
\begin{remark}
We obtain detailed textbook descriptions of
 ${\rm trop}(\SH(2,5,0))$  from those for
 ${\rm trop}({\rm Gr}(3,6))$ 
in \cite[Sections 4.4 and 5.4]{MStrop}.
This was pointed out in \cite[Section 6]{BDG}.
Similarly, ${\rm trop}_+(\SH(3,7,0))$
arises from ${\rm trop}_+({\rm Gr}(4,8))$.
The latter fan was studied in \cite[Section~2]{DFGK}.
\end{remark}

We finally turn to the tropical Mandelstam variety.
Recall that ${\rm M}(k,n,r)$ is
 the image of ${\rm SH}(k,n,r)$ under
the Hadamard product map $s$ in (\ref{eq:hadamardmap}).
The tropicalization of this map,
\begin{equation}
\label{eq:summap}
{\rm trop}(s) \,\,: \,\,
\RR^{\binom{n}{k}}/ \RR {\bf 1} \,\times \,
\RR^{\binom{n}{k}}/ \RR {\bf 1} \,\,\rightarrow \,\,
\RR^{\binom{n}{k}}/ \RR {\bf 1} , 
\end{equation}
 computes the sum of two tropical  Pl\"ucker vectors,
modulo global tropical scaling.
It follows from \cite[Theorem 5.5.1]{MStrop} that the 
Hadamard product map $s$ commutes with tropicalization. 
The following tropical constructions are thus obtained directly from
their classical analogues.

\begin{corollary} The tropical Mandelstam variety is the
image of  the tropical spinor-helicity variety under the 
sum map (\ref{eq:summap}).
Namely, for all values of the parameters $k,n,r$, we~have
 \begin{equation}
 \label{eq:tropM}
 {\rm trop}({\rm M}(k,n,r)) \,\, = \,\,
{\rm trop}(s) \bigl( \,{\rm trop}({\rm SH}(k,n,r)) \,\bigr). 
\end{equation}
For special parameter values,
the tropical Mandelstam variety
is the image of a tropical Grassmannian or a tropical flag variety  under the sum map. 
For instance, the largest
tropical Mandelstam variety ${\rm M}(k,n,k)$ is simply
the Minkowski sum of the Grassmannian with itself:
\begin{equation}
 \label{eq:tropM2}
 {\rm trop}({\rm M}(k,n,k)) \,\, = \,\,
 {\rm trop}({\rm Gr}(k,n))\,+\,
 {\rm trop}({\rm Gr}(k,n)). 
 \end{equation}
\end{corollary}
It would be an interesting future project
to study the combinatorics of these polyhedral fans.

\begin{example}[$k=2,n=5$]
Equation (\ref{eq:tropM}) describes a $2$-to-$1$ map
from ${\rm trop}({\rm Gr}(3,6))$ onto ${\rm trop}({\rm M}(2,5,0))$.
Starting from the census in \cite[Example 4.4.10]{MStrop},
we can examine this map on every cone of ${\rm trop}({\rm Gr}(3,6))$.
This fan has
$65 = 20+15+30$ rays, grouped into types E, F and G.
These rays map to the $40 = 10+15+15$ rays of ${\rm trop}({\rm M}(2,5,0)) = 
{\rm trop}(R_{10})$.
\end{example}

\section{The Scattering Correspondence}
\label{sec6}

The {\em scattering potential} in the
CEGM model \cite{CEGM} is
$$
L_s \,\, = \,\,\sum_{I \in \binom{[n]}{k}} s_I \cdot {\rm log}( p_I ) .
$$
We assume that  $s = (s_I)$ is a fixed point
in the Mandelstam variety ${\rm M}(k,n,r)$ where
$r \leq k-2$. The unknowns $p = (p_I)$ are
the Pl\"ucker coordinates of the open
Grassmannian ${\rm Gr}(k,n)^o$, 
which is defined by $p_I \not= 0$ for all $I \in \binom{[n]}{k}$.
The momentum conservation relations on
${\rm M}(k,n,k-2)$ ensure that
the scattering potential is a well-defined on the quotient space
\begin{equation}
\label{eq:quotientspace}
 X(k,n) \,\, = \,\, {\rm Gr}^{\circ}(k,n) / (\CC^*)^n . 
 \end{equation}
This is a very affine variety of dimension $(n-k-1)(k-1)$; see \cite{ABF, Lam2}.
It is the moduli space of configurations of $n$ labeled points in
linearly general position in the projective space $\PP^{k-1}$.

The scattering potential $L_s$ serves as 
log-likelihood function in algebraic statistics \cite{ST}.
In both statistics and physics, one cares about the critical points of
$L_s$. These are defined~by
\begin{equation}
\label{eq:scatteqn}
\nabla_p \, L_s \,\,= \,\, 0 .
\end{equation}
We now let both $s$ and $p$ vary, and we consider
all solution pairs $(s,p)$ to the system of equations in (\ref{eq:scatteqn}).
The pairs $(s,p)$ satisfying  (\ref{eq:scatteqn})
 are the points of the {\em scattering correspondence}
$$ {\rm C}(k,n,r) \,\, \subset \,\, {\rm M}(k,n,r) \,\times\, {\rm X}(k,n). $$
The aim of this section is to initiate the mathematical study of
this variety. We also consider the 
{\em lifted scattering correspondence}, where ${\rm M}(k,n,r)$ is 
replaced by the spinor-helicity variety:
$$ \tilde {\rm C}(k,n,r) \,\, \subset \,\, {\rm SH}(k,n,r) \,\times\, {\rm X}(k,n). $$

\begin{example}[$k=3,n=6$]
We represent six points in $\PP^2$ by the columns of a matrix
\begin{equation} 
\label{eq:P}
P \quad = \quad \begin{bmatrix}
                        1 &   0  &  0 &   1 &   1 &   1 \\
                        0  &  1  &  0  &  1  &  x  &  y \\
                        0  & 0  &  1 &   1 &   z &   w 
\end{bmatrix}
\end{equation}
Hence $x,y,z,w$ are coordinates on the moduli space ${\rm X}(3,6)$.
The scattering potential equals
\[ \resizebox{\textwidth}{!}{$\begin{matrix} L_s & \!\! =\! \! &\,\,\, s_{125}\cdot  {\rm log}(z) 
+ s_{126}\cdot {\rm log}(w) 
+ s_{135} \cdot {\rm log}(-x) + s_{136} \cdot {\rm log}(-y) + s_{145} \cdot {\rm log}(z-x)
     + s_{146}\cdot {\rm log}(w - y)  \smallskip \\
 & &   \!\!\!\!  \! + \,\,s_{156} \cdot {\rm log}(w x - y z) + s_{245} \cdot {\rm log}(1-z ) 
 + s_{246}\cdot  {\rm log}(1-w)
     + s_{256}\cdot  {\rm log}(z-w) + s_{345}\cdot {\rm log}(x - 1)\smallskip \\ & & 
      + \,\,s_{346}\cdot {\rm log}(y - 1) \,+\, s_{356}\cdot {\rm log}(y-x)
           \,+\, s_{456} \cdot {\rm log}(w x - y z - w - x + y + z).
\end{matrix}$}
\]
The scattering equations are given by the partial derivatives of the scattering potential $L_s$:
\begin{equation}
\label{eq:scatter36}
 \begin{matrix} s_{126}\, \frac{1}{w} + s_{146}\, \frac{1}{w-y} +
s_{156} \, \frac{x}{ wx-yz} - s_{246}\, \frac{1}{1-w} - s_{256}\,\frac{1}{z-w} 
+ s_{456}\, \frac{x-1}{wx - yz -w-x+y+z} & = & 0,\smallskip \\
s_{135}\,\frac{1}{x} - s_{145} \,\frac{1}{z-x} + s_{156}\,\frac{w}{wx-yz} + s_{345} \, \frac{1}{x-1}
- s_{356}\, \frac{1}{y-x} + s_{456}\, \frac{w-1}{wx - yz -w-x+y+z} & = & 0, \smallskip \\
s_{136}\,\frac{1}{y} - s_{146} \,\frac{1}{w-y} - s_{156}\,\frac{z}{wx-yz} + s_{346} \,\frac{1}{y-1}
+ s_{356} \frac{1}{y-x} + s_{456}\, \frac{1-z}{wx - yz -w-x+y+z} & = & 0,\smallskip \\ \,
s_{125}\,\frac{1}{z} + s_{145}\,\frac{1}{z-x}\, - \,s_{156}\, \frac{y}{wx-yz} - s_{245} \,\frac{1}{1-z}
+ s_{256} \,\frac{1}{z-w} + s_{456}\, \frac{1-y}{wx - yz -w-x+y+z} & = & 0.
\end{matrix} 
\end{equation}

If the $s_{ijk}$ are general solutions to the
linear equations in (\ref{eq:MCR361}) then
(\ref{eq:scatter36}) has $26$ complex solutions in ${\rm X}(3,6)$.
We are interested in the case when $s$ lies in ${\rm M}(3,6,1)$, or when we  lift to
${\rm SH}(3,6,1)$ by substituting $s_{ijk} = \langle ijk \rangle [ijk]$.
We obtain $26$-to-$1$ maps from the two scattering
correspondences ${\rm C}(3,6,1)$ or $\tilde {\rm C}(3,6,1)$
onto their kinematic spaces ${\rm M}(3,6,1)$ or ${\rm SH}(3,6,1)$.
      \end{example}

We next examine the case of primary interest in physics, namely $k=2$.
The Mandelstam variety ${\rm M}(2,n,r)$ was characterized in
Section \ref{sec4}.
Here, $X(2,n)$ is the moduli space $\mathcal{M}_{0,n}$
of $n$ distinct labeled points $x_1,x_2,\ldots,x_n$ in $\PP^1 = \CC \cup \{\infty\}$.
The scattering potential equals
$$ L_s \,\,\, = \, \sum_{1 \leq i < j \leq n} \! \! s_{ij} \cdot {\rm log}( x_i - x_j). $$
The system of scattering equations $\nabla_x L_s = 0$ can be written explicitly as follows:
\begin{equation}
\label{eq:scat2b} \sum_{j=1}^n \frac{s_{ij}}{ x_i - x_j} \,\, = \,\, 0 
\qquad {\rm for} \,\, i = 1,2,\ldots,n . 
\end{equation}

Let $z$ be a new unknown and consider the following rational function in $z$
of degree $-2$:
\begin{equation}
\label{eq:scat2c}
T(z) \,\,\, = \,\sum_{1 \leq i < j \leq n} \frac{s_{ij}}{ (z-x_i)(z-x_j)}.
\end{equation}

\begin{proposition}
The rational function $T(z)$ is identically zero if and only if
(\ref{eq:scat2b}) holds.
\end{proposition}

\begin{proof}
The residue of $T(z)$ at $z=x_i$ is precisely the left hand side of
the equation in (\ref{eq:scat2b}).
The residues at the $n$ poles are all zero if and only if $T(z)$ is the zero function.
\end{proof}

This following result is known in the particle physics literature due to
work of Witten, Roiban-Spradlin-Volovich and Cachazo-He-Yuan \cite{CHY13}.
We learned it from  the recent 
lectures by Thomas Lam \cite[Section 4.4]{Lam2}.
See the discussion of ``sectors'' in \cite[Introduction]{Lam2}.

\begin{theorem} \label{thm:LSC}
The lifted scattering correspondence  $\tilde {\rm C}(2,n,0)$
has $n-3$ irreducible components
$ \tilde{\rm C}_2, \tilde{\rm C}_3, \ldots, \tilde{\rm C}_{n-2}$. Each of them has
 the same dimension as ${\rm SH}(2,n,0)$.
The irreducible component $\tilde{\rm C}_\ell$ parametrizes all $3$-step flags
$V \subseteq U \subseteq W^\perp$
where $(V,W) $ is a point in $ {\rm SH}(2,n,0)$ 
and $U $ is the row span of an
$\ell \times n$ Vandermonde matrix  $(x_j^i)_{i=0,\ldots,\ell-1; j=1,\ldots,n}$.
The map from $\tilde{\rm C}_{\ell}$ to ${\rm SH}(2,n,0)$ is 
finite-to-one: its degree  is the Eulerian number $A(n-3,\ell-2)$.
\end{theorem}

\begin{remark}
The maximum likelihood degree of the moduli space $\mathcal{M}_{0,n}$ equals $(n-3)!$.
In other words, the equations (\ref{eq:scat2b}) have $(n-3)!$ solutions, provided
$\sum_{j=1}^n s_{ij} = 0$ for all $i$. See e.g.~\cite[Proposition 1]{ST}. 
Theorem~\ref{thm:LSC}  is a geometric realization of the combinatorial identity
\begin{equation}
\label{eq:factorialsum} (n-3)! \,\, = \,\, A(n-3,0) + A(n-3,1) +\, \cdots\, + A(n-3,n-4). 
\end{equation}
Note that the Eulerian numbers can be defined by
 $ A(2,0) = A(2,1) = 1 $ and the recursions
 \begin{equation}
 \label{eq:Arecursion}
 A(n-3,\ell-2)\,\, =\,\, (\ell-1) \cdot A(n-4,\ell-2) \,+\, (n-\ell-1) \cdot A(n-4,\ell-3) 
 \quad {\rm for} \,\, n \geq 3. 
 \end{equation}
 \end{remark}
\begin{proof}[Proof of Theorem \ref{thm:LSC}]
A proof for the first part of the statement concerning the irreducible components of $\tilde{\rm C}(2,n,0)$ was given in \cite[Proposition 4.6 and Section 4.5]{Lam2}. The argument
for the second part was outlined in \cite[Section 1.1]{CHY13}. It is based  on
the degeneration technique known in physics as {\em soft limits}.
We now present details from the algebraic perspective of \cite{ABF}.
See also \cite[Corollary 4.8]{ARS} for a recent proof using intersection theory.

The variety $\tilde {\rm C}_\ell$ is the image of
the following map to the lifted scattering correspondence:
\begin{equation}
\label{eq:LSCpara}
 (\CC[z]_{\leq \ell-1})^2 \times (\CC[z]_{\leq n-\ell-1})^2 \times \mathcal{M}_{0,n}
\,\rightarrow \,\tilde {\rm C}(2,n,0) \,,\,\,
(\tau, \tilde \tau, x) \,\mapsto\, (\lambda, \tilde \lambda, x) . 
\end{equation}
Here $x = (x_1,\ldots,x_n) \in \CC^n $ represents a point in $\mathcal{M}_{0,n}$,
the column vector $\tau$ resp.~$\tilde \tau$ consists of two polynomials
in one variable $z$ of degrees $\ell-1$ resp.~$n-\ell-1$,
the $i$th column of  $\lambda$
equals $\,\tau(x_i)\,$, and the $i$th column of $\,\tilde \lambda\,$ is
$\,\tilde{\tau}(x_i) \prod_{j \not= i} (x_i-x_j)^{-1} $.
Note that the row spaces $V$ and $W$ of these $2 \times n$ matrices satisfy
$V \subseteq U \subseteq W^\perp$, where $U$ is the row space of $(x^i_j)$.

Write $D(n,\ell)$ for the degree of the finite-to-one map
$\,\tilde {\rm C}_\ell \rightarrow {\rm SH}(2,n,0),\, (\lambda,\tilde \lambda,x) \mapsto
(\lambda,\tilde \lambda)$.
In words, $D(n,\ell)$ is the number of solutions to the
scattering equations in the $\ell$th component of 
$\tilde{\rm C}(2,n,0)$. Direct computations reveal
$D(4,2) = D(5,2) = D(5,3) = 1$. We claim that
 \begin{equation}
 \label{eq:Drecursion}
 D(n,\ell) \,\,= \,\,(\ell-1)\cdot D(n-1,\ell) \,+\, (n-\ell-1)\cdot D(n-1,\ell-1)
 \quad {\rm for} \,\, n \geq 3.
 \end{equation}
This claim implies Theorem \ref{thm:LSC}, because
the recursion in (\ref{eq:Drecursion}) matches the recursion in
(\ref{eq:Arecursion}).

It remains to prove (\ref{eq:Drecursion}).
This is done using the technique of soft limits.
We drive the $n$th particle to zero, by replacing
the Mandelstam coordinate
$s_{in}$ with $ \epsilon \cdot s_{in} = \epsilon \cdot \langle in \rangle [in]$ for $1 \leq i \leq n-1$.
This degeneration is compatible with the parametrization
(\ref{eq:LSCpara}). As $\epsilon$ tends to zero,
either $\tau(x_n)  \rightarrow 0$ or $\tilde \tau(x_n) \rightarrow 0$.
This yields one equation in one unknown $x_n$ of degree $\ell-1$ resp.~$n-\ell-1$.
The other unknowns $x_1,\ldots,x_{n-1}$ satisfy the scattering equations
on the components of $\tilde {\rm C}(2,n-1,0)$ that are
indexed by $\ell$ and $\ell-1$ respectively.
So, for $\epsilon$ near $0$, the size of the fibers of
$\,\tilde {\rm C}_\ell \rightarrow {\rm SH}(2,n,0) \,$
is the right hand side of (\ref{eq:Drecursion}). Finally, we note that for generic kinematics $(s_{ij}) \in {\rm M}(2,n,0)$, the zero-dimensional polynomial system at hand is reduced, and the last part of the argument rests on applying the Implicit Function Theorem.
\end{proof}

 The irreducible components $\tilde {\rm C}_i$
 are referred to as ``sectors'' in the physics literature.
  Section 5.1 in \cite{CEGM} starts with the sentence
 {\em
``In the k = 2 case it is well known that solutions of the scattering equations split into $n-3$ sectors''}.
Our proof was written for mathematicians.

\smallskip

Theorem \ref{thm:LSC}
was stated for the  lifted scattering correspondence $\tilde {\rm C}(2,n,0)$.
From our perspective, it is more natural to focus on
the scattering correspondence ${\rm C}(2,n,0)$
because the parameters in $L_s$ are the Mandelstam invariants.
The scattering correspondence ${\rm C}(2,n,0)$ lives over
the Mandelstam variety ${\rm M}(2,n,0)$, whose
prime ideal we presented in Theorem~\ref{thm:prime2n}.

\begin{corollary} \label{cor:SC}
The scattering correspondence  $ {\rm C}(2,n,0)$
has $\lceil \frac{n-3}{2} \rceil$ irreducible components.
The varieties ${\rm C}_\ell$ and ${\rm C}_{n-\ell}$
in Theorem \ref{thm:LSC}  are identified by the
map ${\rm SH}(2,n,0) \rightarrow {\rm M}(2,n,0)$.
\end{corollary}

\begin{proof}
The Hadamard product map $s$ gives rise to a
map from the lifted scattering correspondence
in  ${\rm SH}(2,n,0) \times \mathcal{M}_{0,n}$ 
onto the scattering correspondence in $ {\rm M}(2,n,0) \times \mathcal{M}_{0,n}$.
This map is a covering of degree two. 
The fibers are comprised of solutions for $\ell$ and for~$n - \ell$.
These are distinct, unless $ \ell = n/2$, where
the solution is fixed under this involution.
 \end{proof}
 
We now turn to the case $k=3, r=1$, which was 
discussed by Cachazo et al.~in \cite[Section 5.1]{CEGM}.
We revisit their results from an algebraic perspective,
and we report on the identification of irreducible components
with the help of {\tt HomotopyContinuation.jl} \cite{BT}.

The moduli space $X(3,n)$ is a quotient of 
the Grassmannian ${\rm Gr}(3,n)$,
by  (\ref{eq:quotientspace}). Hence
there are {\em two tautological maps}
from ${\rm SH}(3,n,1)$ to $X(3,n)$.
These rational maps take 
the pair $(V,W)$ in (\ref{eq:rephrasing})
 to the images of $V$ 
and $W$ in $X(3,n)$ respectively.
Furthermore, we consider the Veronese map $\nu$ from
$ \mathcal{M}_{0,n}$ into $X(3,n)$
which take $n$ points in $\PP^1$ 
to $n$ points on a conic in $\PP^2$.
Algebraically, the map $\nu$ takes
a $2 \times n$ matrix with $i$th column $(u_i,v_i)^\top$
to the $3 \times n$ matrix with $i$th column $(u_i^2, u_i v_i, v_i^2)^\top$.
Its Pl\"ucker coordinates satisfy
$\langle ijk \rangle \,= \, \langle ij \rangle  \langle ik \rangle \langle jk \rangle$.

We now compose the rational map in (\ref{eq:k-2map}) with 
the two tautological maps for $k=2$, followed by the
Veronese map $\nu$. In this manner, we obtain 
{\em two tautological Veronese~maps}
$$ {\rm SH}(3,n,1) \,\,\dashrightarrow\,\,
{\rm SH}(2,n,0)\, \,\dashrightarrow \,\, X(2,n) \,=\, \mathcal{M}_{0,n}\,
\,\stackrel{\nu}{\hookrightarrow}\,\, X(3,n).$$
Explicitly, these two rational maps are
$\,(V,W) \mapsto \nu( V \cap W^\perp) \,$ and
$\,(V,W) \mapsto \nu( V^\perp \cap W) $.
We paraphrase the construction in \cite[Section 5.1]{CEGM} 
in terms of the scattering correspondence.

\begin{theorem}
The lifted scattering correspondence $\tilde {\rm C}(3,n,1)$
contains four irreducible components which map 
birationally onto the moduli space $X(3,n)$.
These are given by the
two tautological maps and the two tautological Veronese maps.
In other words, for general kinematic data $(V,W) \in \SH(3,n,1)$,
with Mandelstam invariants $s \in {\rm M}(3,n,1)$,
the configurations $V,\,W,\, \nu( V \cap W^\perp) ,\, \nu( V^\perp \cap W) $ in
$ X(3,n)$ are solutions to the scattering equations~(\ref{eq:scatteqn}).
\end{theorem}

We illustrate this result by solving the equations
(\ref{eq:scatter36}) for a random numerical instance.

\begin{example}[$k=3,n=6$]
Let $(V,W) \in {\rm SH}(3,6,1)$ be the  point with parameter matrix 
\[
{\bf x} \,\,\, = \,\,\, \begin{footnotesize}
\begin{bmatrix}
                           4 &   0 &   7 &   4  &  9  &  1 \\
                           1 &   3 &   7 &   2  &  8  &  9 \\
                            1 &   7 &   9 &   8 &   0 &   5 \\
                            6  &  6  &  2  &  2  &  4  &  2
\end{bmatrix}. \end{footnotesize}
\]
Thus $V$ is the span of the first three rows and
$W$ is the kernel of the last three rows.
The Mandelstam invariants are computed by
composing the map $\phi_{3,6,1}$ with the
Hadamard map:
$$ \begin{scriptsize} \begin{matrix}
s_{123} = 12000 & s_{124} = 6720 & s_{125} = 8272 & s_{126} = -31584 & s_{134} = -37760 &
    s_{135} = 54784 & s_{136} = -35728 \\ s_{145} = -38080 & s_{146} = 92208& s_{156} = -30832 & 
    s_{234} = -37920 & s_{235} = 68288 & s_{236} = -108016 & s_{245} = -82896 \\ s_{246} = 82416 & 
    s_{256} = 82720 & s_{345} = 46592 & s_{346} = 57664 & s_{356} = -19904 & s_{456} = -88944.
\end{matrix} \end{scriptsize}
$$
The scattering equations $ (\ref{eq:scatter36})$ have $26$ solutions
$(x,y,z,w) \in \CC^4$. Among these, we find:
\begin{center}
\begin{tabular}{|c|c|c|c|} 
 \hline
 $V$ & $W$ & $\nu(V \cap W^\perp)$ & $\nu(V^\perp \cap W)$ \\ [0.5ex] 
 \hline
 \hline 
 $\bigl(\frac{8453}{5723}, \frac{6083}{9263}, \frac{3713}{1358}, \frac{3713}{2198} \bigr)$ 
 & $\bigl(\frac{6}{11}, \frac{87}{172}, -\frac{1}{4}, -\frac{42}{43} \bigr)$ 
 & $\bigl(\frac{5}{21}, \frac{5}{36}, \frac{19}{27}, \frac{38}{69} \bigr)$ 
 & $\bigl(\frac{6588}{14911}, \frac{20988}{9139}, -\frac{8601}{125060}, \frac{17437}{138528} \bigr)$\\ [1ex] 
 \hline
\end{tabular}
\end{center}
\medskip
\noindent
Substituting these $(x,y,z,w)$ into (\ref{eq:P}), we
obtain four $3 \times 6$ matrices $P$ with rational entries.
These matrices represent the configurations
$V,\,W,\, \nu( V \cap W^\perp) ,\, \nu( V^\perp \cap W)$, where $V$ is the linear span of the first three rows of $X$ and $W^\perp$ is the span of the last three rows of $X$.
\end{example}

\begin{conjecture} \label{conj:end}
The lifted scattering correspondence 
$\tilde {\rm C}(3,n,1)$ decomposes into five irreducible components,
i.e.~there is only one component whose map
 onto $X(3,n)$ is not birational.
Hence, the scattering correspondence  ${\rm C}(3,n,1)$
 has three irreducible components.
 The four birational components become two components modulo
 the involution in Remark~\ref{rmk:involution}.
 \end{conjecture}

The passage from the first sentence to the second sentence
in Conjecture \ref{conj:end} mirrors the passage from
Theorem \ref{thm:LSC}
to Corollary \ref{cor:SC}.
We verified our conjecture in some small cases.
In particular, we showed that
Conjecture \ref{conj:end} is true for $n=6,7,8$.
The verification in these three cases is a computation with 
the software {\tt HomotopyContinuation.jl} \cite{BT}.
Recall from \cite{ABF} that the degree of
the map $\tilde {\rm C}(3,n,1) \rightarrow \SH(3,n,1)$  equals
$26$, $1272$, $188112$ for $n=6,7,8$.
We ran numerical irreducible decomposition,
based on monodromy loops, on the defining equations
of the lifted scattering correspondence.
We found five irreducible components over $X(3,n)$.
For instance, for $n=7$, the five components have
degrees $1,1,1,1$ and $1268$.

\section*{Acknowledgments}
YEM is supported by Deutsche Forschungsgemeinschaft (DFG, German Research Foundation) SFB-TRR 195 ``Symbolic Tools in Mathematics and their Application''.
AP and BS are supported by the European Union (ERC, UNIVERSE+, 101118787).
We are grateful to Simon Telen for helping us with numerical computations for Section \ref{sec6}.

\bibliographystyle{siamplain}

\begin{thebibliography}{10}
    
    \bibitem{ABF} D.~Agostini, T.~Brysiewicz, C.~Fevola, L.~K\"uhne, B.~Sturmfels and S.~Telen:
      {\em Likelihood degenerations}, Adv.~Math.~{\bf 414} (2023) 108863.

    \bibitem{ARS}
    D.~Agostini, L.~Ramesh and D.~Shen:
    Points on rational normal curves and the ABCT variety, 
    {Le Matematiche} {\bf 80} (2025) 103-122.
    
    \bibitem{AKW}
    F.~Ardila, C.~Klivans and L.~Williams:
    {\em  The positive Bergman complex of an oriented matroid},
    Eur. J. Comb.~{\bf 27} (2006) 577--591.
    
    
    \bibitem{BHPZ}
    S.~Badger, J.~Henn, J.~Plefka and S.~Zoia:
    {\em Scattering Amplitudes in Quantum Field Theory},
    Lecture Notes in Phys. {\bf 1021}, Springer,  2024.
    
    \bibitem{BPT}
    B.~Betti, M.~Panizzut and S.~Telen:
    {\em Solving equations using Khovanskii bases},
    J. Symbolic Comput. 126 (2025), Paper No. 102340.
    
    \bibitem{Bloch_Karp}
    A.~Bloch and S.~Karp:
    {\em On two notions of total positivity for partial flag varieties},
    Adv. Math. 414 (2023) 108855.
    
    \bibitem{BC}
    C.~Bocci and E.~Carlini:
    {\em  Hadamard Products of Projective Varieties},
    Birkh\"auser, Cham, 2024.
     
    \bibitem{Bor}
    J.~Boretsky: {\em Positive tropical flags and the positive tropical Dressian},
    S\'emin. Lothar. Comb. 86B (2022) Article 86.
     
    \bibitem{BDG}
    L.~Bossinger, J.~Drummond and R.~Glew:
    {\em Adjacency for scattering amplitudes from the Gr\"obner fan},
    J.~High Energy Phys.~2023, no. 11, 2.
    
    \bibitem{BLS}
    M.~Brandenburg, G.~Loho and R.~Sinn:
    {\em Tropical positivity and determinantal varieties},
    Algebr. Comb. {\bf 6} (2023) 999--1040.
    
    \bibitem{BEZ}
    M.~Brandt, C.~Eur and L.~Zhang: {\em Tropical flag varieties},
    Adv.~ Math.~{\bf 384} (2021) 107695.
    
    \bibitem{BT} P.~Breiding and S.~Timme:
    {\em HomotopyContinuation.jl: A package for homotopy continuation in Julia},
    Mathematical Software -- ICMS 2018, 458--465,
     Springer, 2018.
     
    \bibitem{CEGM} F.~Cachazo, N.~Early, A.~Guevara and S.~Mizera: 
    {\em Scattering equations: from projective spaces to tropical Grassmannians},
    J.~High Energy~Phys.~(2019), no. 6, 39.
    				
    \bibitem{CE} F.~Cachazo and N.~Early:
    {\em Biadjoint scalars and associahedra from residues of generalized amplitudes},
    J.~High Energy~Phys.~(2023), no. 10, 15.
    
    \bibitem{CHY13}{F.~Cachazo, S.~He and E.~Yuan: {\em Scattering in three dimensions from rational maps}, J. High Energy Phys. (2013) no. 141.}
    
    \bibitem{DFRS}
    K.~Devriendt, H.~Friedman, B.~Reinke and B.~Sturmfels:
    {\em The two lives of the Grassmannian},
    Acta Univ. Sapientiae, Math. (2025).
    
    \bibitem{DFGK}
    J.~Drummond, J.~Foster, \"O.~G\"urdo\c{g}an and C.~Kalousios:
    {\em Algebraic singularities of scattering amplitudes from tropical geometry},
    J.~High Energy~Phys.~(2021), no 4, 2.
    
    \bibitem{EH} H.~Elvang and Y.~Huang: {\em Scattering Amplitudes 
    in Gauge Theory and Gravity}, Cambridge University Press, 2015.
    
    \bibitem{GM16}
    C.~Godsil  and K.~Meagher: {\em Erd\H{o}s-{K}o-{R}ado Theorems: Algebraic Approaches},
    Cambridge University Press, 2016.
    
    \bibitem{M2} D.~Grayson  and M.~Stillman:  Macaulay2, a software system for research in algebraic geometry, available at
    {\tt https://macaulay2.com}.
    
    \bibitem{LG} J.~Huh and B.~Sturmfels:  {\em Likelihood geometry},
    in Combinatorial Algebraic Geometry (eds. Aldo Conca et al.), Lecture Notes in Math. {\bf 2108},
      Springer Verlag, (2014) 63--117.
    
    \bibitem{Lam1} T.~Lam: {\em An invitation to positive geometries}, {Open problems in algebraic combinatorics}, Proc. Sympos. Pure Math., 110, American Mathematical Society 2024.

    \bibitem{Lam2} T.~Lam: {\em Moduli spaces in positive geometry},
    {Le Matematiche} {\bf 80} (2025) 17-101
        
    \bibitem{MStrop}  D.~Maclagan and B.~Sturmfels:  {\em Introduction to Tropical Geometry},
     Grad. Stud. Math., Vol 161, American Mathematical Society, 2015.
    
    \bibitem{MS} M.~Micha\l ek and B.~Sturmfels:  {\em Invitation to Nonlinear Algebra},
     Grad. Stud. Math., Vol 211, American Mathematical Society, 2021.
    
    \bibitem{CCA}
    E.~Miller and B.~Sturmfels: {\em Combinatorial Commutative Algebra},
    Grad. Texts in Math., Vol 227, Springer-Verlag, New York, 2005.
    
    \bibitem{Ola} J.~Olarte: {\em Positivity for partial tropical
    flag varieties}, {\tt arXiv:2302.10171}.
    \bibitem{RSS} S.~Rajan, S.~Sverrisd{\'o}ttir and B.~Sturmfels: {\em Kinematic varieties for massless particles}, Le Mathematiche  {\bf 80} (2025) 53–471.
    
    \bibitem{SW}
    D.~Speyer and L.~Williams: {\em  The tropical totally positive Grassmannians},
    J. Algebr.~Comb.~{\bf 22} (2005) 189--210.
    
    \bibitem{AIT} B.~Sturmfels:  {\em Algorithms in Invariant Theory}, Springer-Verlag, Vienna, 1993.
    
    \bibitem{ST} B.~Sturmfels and S.~Telen: {\em Likelihood equations and scattering amplitudes}, Algebr.~Stat.~{\bf 12} (2021) 167--186.
    
    \bibitem{YY} J.~Yu and D.~Yuster:
    {\em Representing tropical linear spaces by circuits},
    Formal Power Series and Algebraic Combinatorics (FPSAC 2007),
    {\tt arXiv:math/0611579}.
\end{thebibliography}

\end{document}